\documentclass{amsart}
\usepackage{amsmath}
\usepackage{amssymb}
\usepackage{latexsym}
\usepackage{amsfonts}
\usepackage[dvips]{graphicx}
\usepackage{color}
\usepackage{graphicx,epstopdf}

\newtheorem{thm}{Theorem}

\newtheorem{corollary}{Corollary}

\usepackage{float}
\usepackage[section]{placeins}

\usepackage{color}
\usepackage{graphicx}
\usepackage{epstopdf}

\begin{document}

\title[In search for a perfect shape of polyhedra: Buffon transformation]{In search for a perfect shape of polyhedra: Buffon transformation}
\author{V. Schreiber}\address{Department of Mathematical Sciences,
Loughborough University, Loughborough LE11 3TU, UK}
\email{V.Schreiber@lboro.ac.uk}
\author{A.P. Veselov}
\address{Department of Mathematical Sciences,
Loughborough University, Loughborough LE11 3TU, UK  and Moscow State University, Moscow 119899, Russia}
\email{A.P.Veselov@lboro.ac.uk}
\author{J.P. Ward}
\address{Department of Mathematical Sciences,
Loughborough University, Loughborough LE11 3TU, UK}
\email{J.P.Ward@lboro.ac.uk}

\maketitle

\begin{abstract}
For an arbitrary polygon consider a new one by joining
the centres of consecutive edges. Iteration of this procedure leads to
a shape which is affine equivalent to a regular polygon.
This regularisation effect is usually ascribed to Count Buffon (1707-1788).
We discuss a natural analogue of this procedure for 3-dimensional polyhedra,
which leads to a new notion of affine $B$-regular polyhedra.
The main result is the proof of existence of star-shaped affine $B$-regular polyhedra 
with prescribed combinatorial structure, under partial symmetry and simpliciality 
assumptions. The proof is based on deep results from spectral graph theory due to Colin de Verdi\`ere and Lov\'asz.
\end{abstract}

\bigskip

\section{Introduction}
According to David Wells \cite{DWells} the following puzzle first
appeared in Edward Riddle's edition (1840) of the {\it Recreations in Mathematics
and Natural Philosophy} of Jacques Ozanam, where it was attributed to Count
Buffon (1707--1788), a French naturalist and the translator
of Newton's Principia. 

Consider an arbitrary polygon. Generate a second
polygon by joining the centres of consecutive edges. Repeat this construction (see Fig. 1).

\begin{figure}
\includegraphics[scale=0.8]{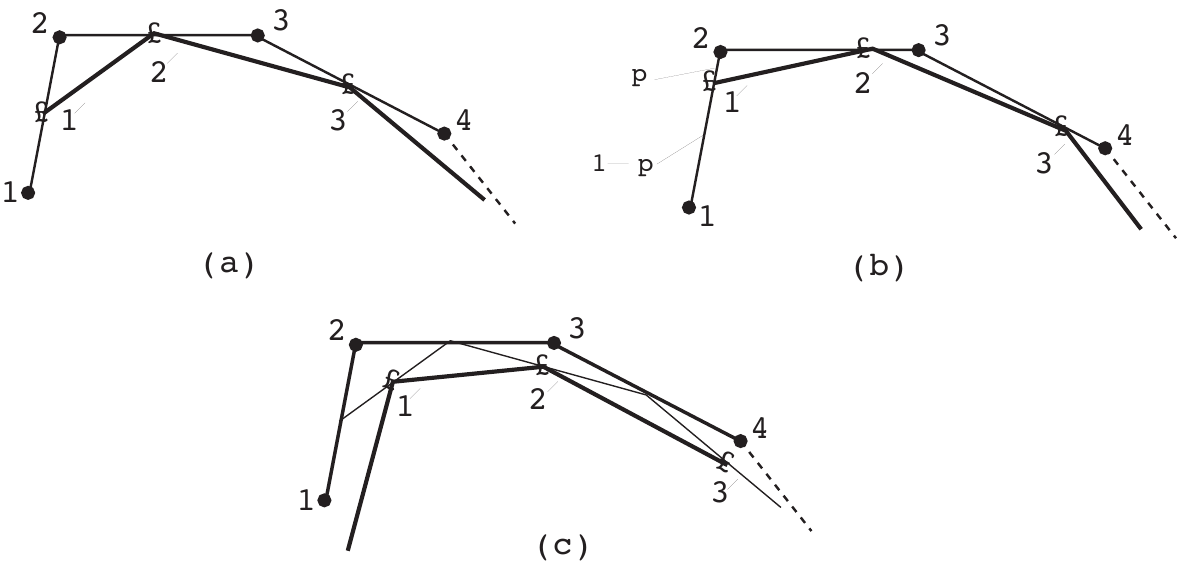}\quad \quad \quad\includegraphics[scale=0.8]{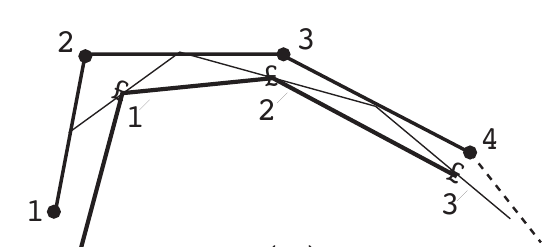}
\caption{} \end{figure}

It is easy to see that the process
converges to a point - the centroid of the original vertices (and therefore
the centroid of the vertices of any polygon in the sequence). 
Buffon observed a remarkable {\it regularization effect} of this procedure:
the limiting shape of the polygon is {\it affine regular.}  
Here a polygon is called affine regular if it is affine
equivalent to a regular polygon. 

In fact such a similar phenomenon was already
observed since Roman times. When creating mosaics Roman craftsmen achieved
more regular pieces by breaking the corner, so effectively using the same procedure \cite{MOS}.
The explanation is based on simple arguments from linear algebra, see e.g. \cite{BGS, W1} and next section.

The situation here is different from the theory of the pentagram map, initiated by R. Schwartz in 1990s and extensively 
studied in recent years, where the dynamics is nonlinear, quasi-periodic and integrable in Arnold-Liouville sense 
(see \cite{KS,OST} and references therein).

In this paper we will study the following natural 3-dimensional version of the Buffon procedure \cite{VW}.
Let $P$ be a simplicial  polyhedron in $\mathbb{R}^{3}$, having all faces triangular.
Define its {\it Buffon transformation} $B(P)$ as the simplicial polyhedron with vertices $B(v)$, 
where for each vertex $v$ of $P$ the new vertex $B(v)$ is defined as the centroid of the centroids of all edges meeting at $v.$  The question is what is the limiting shape of $B^n(P)$ as $n$ goes to infinity.

Unfortunately, the answer in general is disappointing: the limiting shape will be one-dimensional.
Indeed the same arguments from linear algebra show that this shape is determined by the subdominant eigenspace 
of the corresponding operator on the graph $\Gamma(P)$, which is the 1-skeleton of $P$ (see the details below), and this eigenspace generically has dimension 1. This means that in order to have a sensible limiting shape we need to add some assumptions on the initial polyhedron $P.$ 

Let $G$ be one of the symmetry groups $G=T,O,I$ of the Platonic solids: tetrahedron, octahedron/cube, icosahedron/dodecahedron respectively. 
Assume that the combinatorial structure of the initial polyhedron $P$ is $G$-invariant, which means that $G$ faithfully acts on the graph $\Gamma(P).$

Our main result is the following theorem.

\begin{thm}
\label{thm:main} Let $P$ be a simplicial polyhedron in $\mathbb R^3$ with $G$-invariant combinatorial structure. 
Then for a generic $P$ the limiting shape obtained by repeatedly applying Buffon procedure to $P$
is a star-shaped polyhedron $P_B.$ The vertices of $P_B$ are explicitly determined by the subdominant eigenspace of the Buffon operator, 
which in this case has dimension 3.
\end{thm}

The proof is based on the deep results from the spectral theory on graphs due to Colin de Verdi\`ere \cite{CdV} and in particular due to Lov\'asz et al \cite{Holst97thecolin, L1, L2}, who studied the eigenspace realisations of polyhedral graphs. Both assumptions of the theorem, namely simpliciality and platonic symmetry, are essential. 

Recall that the polyhedron $P$ is called {\it star-shaped} (not to be mixed with star polyhedra like Kepler-Poinsot) if there is a point inside it from which one can see the whole boundary of $P$, or equivalently, the central projection gives a homeomorphism of the boundary of $P$ onto a sphere. The precise meaning of the term "generic" will be clear from the next section.

Let us call polyhedron $P$ {\it affine $B$-regular} if $B(P)$ is affine equivalent to $P.$ In dimension 2 this is equivalent to affine regularity (see next section). 
Thus the Buffon procedure produces affine $B$-regular version $P_B$ from a generic polyhedron $P$ with the above properties. As far as we know  the notion of the affine regularity for polyhedra with non-regular combinatorial structures was not discussed in the literature before.

For a generic polyhedron $P$ with combinatorial structure of a Platonic solid the corresponding polyhedron $P_B$ is affine regular, which means that it is affine equivalent to the corresponding Platonic solid. For the Archimedean and Catalan solids however, this is no longer true, see the example of pentakis dodecahedron (dodecahedron with pyramids build on its faces) on Fig.2 and in the Appendix.

\begin{figure}\label{F2}
\centerline{ \includegraphics[scale=0.6]{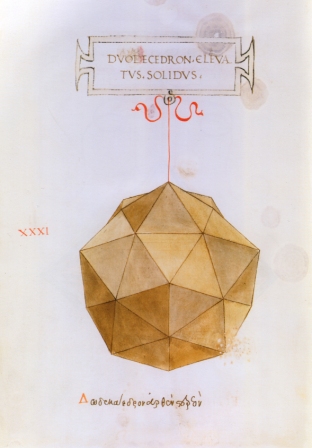}\hspace{20pt}  \includegraphics[scale=0.8]{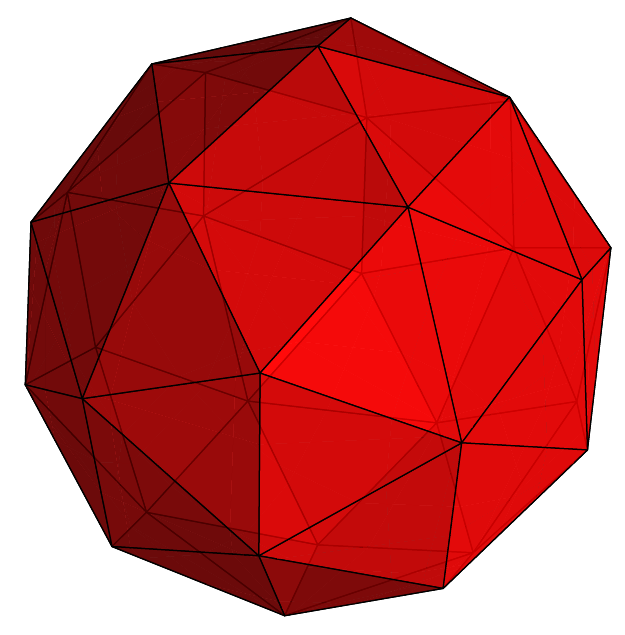}}
\caption{Leonardo da Vinci's drawing of pentakis dodecahedron from Luca Pacioli's book "De divina proportione" and Mathematica image of its Buffon realisation. Leonardo's version is different both for Catalan and Buffon realisations and probably corresponds to the so-called cumulated dodecahedron having all the edges of equal length.} 
\end{figure}

Note that there are plenty of polyhedra $P$ with $G$-invariant combinatorial structures, which can be constructed from the Platonic solids using Conway operations \cite{Con}. In particular, one will have a simplicial polyhedron by applying to any such $P$  the operation, which Conway called $kis$ and denoted by $k$, consisting of building the pyramids on all the faces. Many examples of the corresponding combinatorial types can be found in chemistry and physics literature in relation with the famous Thomson problem, see e.g. \cite{Edmundson}.

For non-simplicial polyhedra the Buffon transformation usually breaks the faces, which in general are not recovering even in the limit (see Fig.13  in Appendix B). 

The platonic symmetry keeps the limiting shape 3-dimensional, preventing collapse to lower dimension.
The dihedral symmetry is not enough: one can check that a polyhedron with prismatic combinatorial structure will collapse to the corresponding affine regular polygon. 

The star-shape property of the limiting shape is probably the strongest we can claim since the convexity may not hold as the example of the triakis tetrahedron shows (see Fig. 11 in the Appendix).

The structure of the paper is as follows. In Section 2 we start with the (well-known) solution of the Buffon puzzle  for polygons to explain the main ideas and relation to linear algebra. Then, in Section 3, we define the Buffon transformation for polyhedra and review the classical Steinitz theorem which gives graph-theoretical characterisation of 1-skeletons of convex polyhedra. In section 4 we introduce the main tools from spectral graph theory: the Colin de Verdi\`ere invariant and null space realisation for polyhedral graphs studied by Lov\'asz et al. In section 5 we use them and representation theory of finite groups to prove our main result. In Appendix A we present the character tables for the polyhedral groups and the corresponding decomposition of the space of functions on the vertices of Platonic solids into irreducible components. In Appendix B we give the spectra of the Buffon operators for some combinatorial types and the corresponding shapes of affine $B$-regular polyhedra.

\section{Buffon transformation for polygons.}

Consider an arbitrary $n$-gon $P$ with vertices described by
the column vector 
\[
r=\left[r_{1},r_{2},\ldots,r_{n}\right], \quad r_i \in \mathbb{R}^{3}
\]
 (and an integer $n\geq3$). Generate a second polygon $P^{'}$
by joining the centres of the consecutive edges of $P$. The corresponding transformation
acts on the vertices of $P$ as follows:
\[
r^{'}_{i}=\frac{1}{2}(r_{i}+r_{i+1}).
\]
In matrix form this can be described as
\[
r^{'}=Br
\]
where
\[
B=\left[\begin{array}{ccccc}
\frac{1}{2} & \frac{1}{2} & 0 & \ldots & 0\\
0 & \frac{1}{2} & \frac{1}{2} & \ldots & 0\\
\vdots & \vdots & \vdots & \vdots & \vdots\\
\frac{1}{2} & 0 & 0 & \ldots & \frac{1}{2}
\end{array}\right]
\]

After $k$
transformations we obtain a polygon with the vertices
\[
r^{k}=B^{k}r.
\]

Following Buffon we claim that for generic initial polygons $P$ 
the limiting shape of the polygons $P^{k}$ as $k$ increases becomes affine regular.
Recall that a polygon is {\it affine regular} if it is affine equivalent to
a regular polygon. 


To prove this we use the following result from Linear Algebra (see e.g. Theorems 5.1.1, 5.1.2 in \cite{Wat}).

\begin{thm} {\bf (Subspace Iteration Theorem)}
Let $A$ be a real matrix and let $Spec(A)=\{\lambda_{1},\lambda_{2},\ldots,\lambda_{n}\}$
be the set of its eigenvalues (in general, complex and with multiplicities) ordered in such a way that
\[
|\lambda_{1}|=|\lambda_{2}|=\ldots=|\lambda_{k}| > |\lambda_{k+1}|\geq\ldots\geq |\lambda_{n}|.
\]
Let $W$ and $W'$ be the dominant and complementary invariant subspaces associated with $\lambda_1, \dots, \lambda_k$ and $\lambda_{k+1}, \dots, \lambda_n$ respectively and $m =\dim W.$
Then for any $m$-dimensional subspace $U \subset \mathbb R^n$ such that
$U\cap W'=\{0\}$ the image of $U$ under the iterations of $A$
\[
A^{n}(U)\underset{n\rightarrow\text{\ensuremath{\infty}}}{\rightarrow}W
\]
tends to the dominant subspace in the Grassmannian $G_{m}(\mathbb R^n)$.
\end{thm}

To apply this to our case first note that 
\[
B=\frac{1}{2}(I+T), 
\]
where the $n\times n$ matrix
\[
T=\left[\begin{array}{ccccc}
0 & 1 & 0 & \ldots & 0\\
0 & 0 & 1 & \ldots & 0\\
\vdots & \vdots & \vdots & \vdots & \vdots\\
1 & 0 & 0 & \ldots & 0
\end{array}\right]
\]
has the property $T^{n}=I$ and the eigenvalues being $n$-th roots of unity.
The spectrum of $B$ is therefore 
\[
Spec(B)=\{\frac{1}{2}+\frac{1}{2}\varepsilon_{j},\:\varepsilon_{j}=e^{\frac{2\pi i}{n}j},j=0,1,\ldots n-1\}.
\]

The eigenvalues of maximum modulus, other than $\lambda_{0}=1$, are
$\lambda_{1}=\frac{1}{2}+\frac{1}{2}e^{\frac{2\pi i}{n}}$ and its
complex conjugate $\lambda_{2}=\frac{1}{2}+\frac{1}{2}e^{-\frac{2\pi i}{n}}=\overline{\lambda_{1}}$.

The dominant subspace $W$ in this case corresponds to $\lambda_0=1$ and is generated by the corresponding eigenvector 
$v_0=(1,1,\ldots,1)$:
$$W=\{(r,r,\dots,r)\}.$$
The previous result can be interpreted that as $n$ increases 
 $B^{n}(P)$ converges to a point. To see the limiting shape we should look at the {\it subdominant invariant} subspace corresponding to $\lambda_1$ and $\lambda_2$.

Geometrically one can do this by assuming that the centroid of the vertices is at the origin ({\it centre of mass condition}). 
This means that we restrict the action of $B$ on the invariant subspace 
$$V_C=\{(r_1, \dots, r_n): r_1+\dots +r_n=0.\}$$
This eliminates the eigenvalue $\lambda_0=1$ and the new dominant subspace $W$
corresponding to $\lambda_{1}=\frac{1}{2}+\frac{1}{2}\varepsilon,\:\lambda_{2}=\overline{\lambda_{1}}$
is precisely the one describing the limiting shape. 
One can easily check that 
\[
W=<\left(\begin{array}{c}
1\\
\varepsilon\\
\varepsilon^{2}\\
.\\
.\\
\varepsilon^{n-1}
\end{array}\right),\left(\begin{array}{c}
1\\
\overline{\varepsilon}\\
\overline{\varepsilon^{2}}\\
.\\
.\\
\overline{\varepsilon^{n-1}}
\end{array}\right)>=\{a\left(\begin{array}{c}
1\\
cos\frac{2\pi}{n}\\
cos\frac{4\pi}{n}\\
.\\
.\\
.
\end{array}\right)+b\left(\begin{array}{c}
0\\
sin\frac{2\pi}{n}\\
sin\frac{4\pi}{n}\\
.\\
.\\
.
\end{array}\right)\}
\]
\medskip{}

Choosing $a$ and $b$ to be orthogonal unit vectors we see that the corresponding vertices form a regular polygon.
In general, the dominant subspace $W$ describes all affine regular polygons.
The other eigenspaces correspond to the affine regular "polygrams".

For example, when $n=5$ we have the eigenvalues
$$\lambda_{1}=\frac{1}{2}+\frac{1}{2}e^{\frac{2\pi i}{n}},\:\lambda_{2}=\overline{\lambda_{1}},\:\lambda_{3}=\frac{1}{2}+\frac{1}{2}e^{\frac{4\pi i}{5}},\:\lambda_{4}=\overline{\lambda_{3}}$$
and the corresponding eigenspaces
$$W=\{a\left(\begin{array}{c}
1\\
cos\frac{2\pi}{5}\\
cos\frac{4\pi}{5}\\
cos\frac{6\pi}{5}\\
cos\frac{8\pi}{5}
\end{array}\right)+b\left(\begin{array}{c}
0\\
sin\frac{2\pi}{5}\\
sin\frac{4\pi}{5}\\
sin\frac{6\pi}{5}\\
sin\frac{8\pi}{5}
\end{array}\right)\},\qquad W^{'}=\{a\left(\begin{array}{c}
1\\
cos\frac{4\pi}{5}\\
cos\frac{8\pi}{5}\\
cos\frac{2\pi}{5}\\
cos\frac{6\pi}{5}
\end{array}\right)+b\left(\begin{array}{c}
0\\
sin\frac{4\pi}{5}\\
sin\frac{8\pi}{5}\\
sin\frac{2\pi}{5}\\
sin\frac{6\pi}{5}
\end{array}\right)\}$$
describing the affine regular pentagons and pentagrams respectively:

\bigskip

\begin{figure} [H]
   \centering
 \includegraphics[scale=0.34]{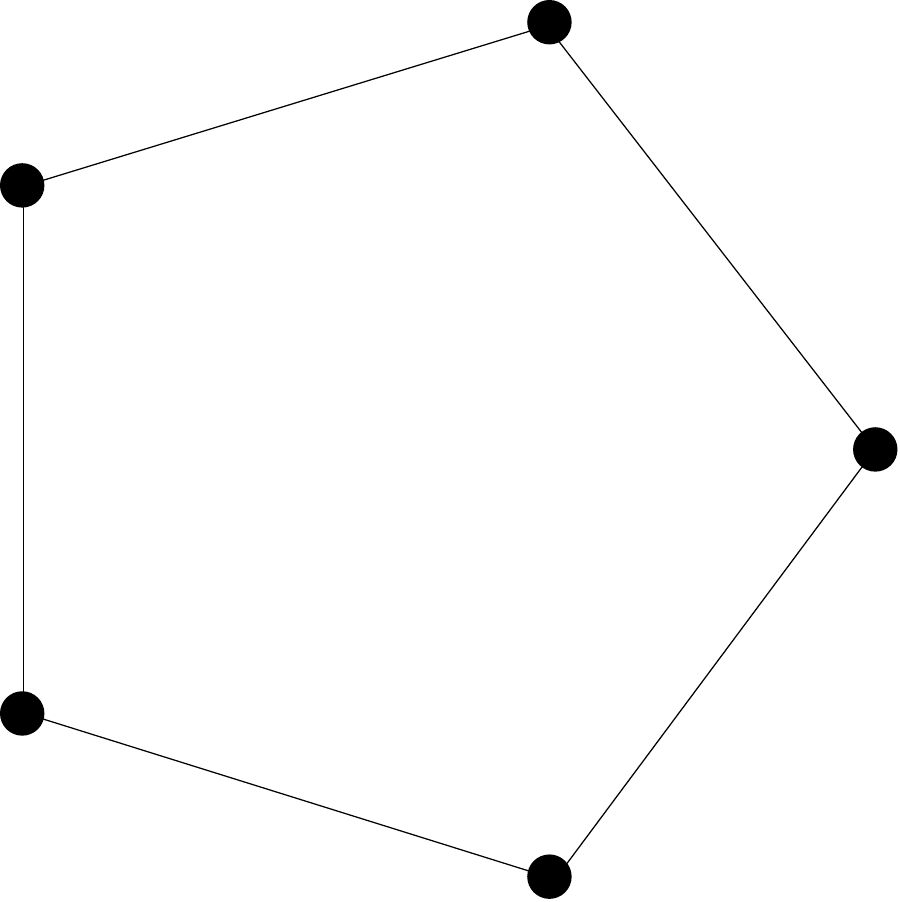}\hspace{2.5cm}\includegraphics[scale=0.4]{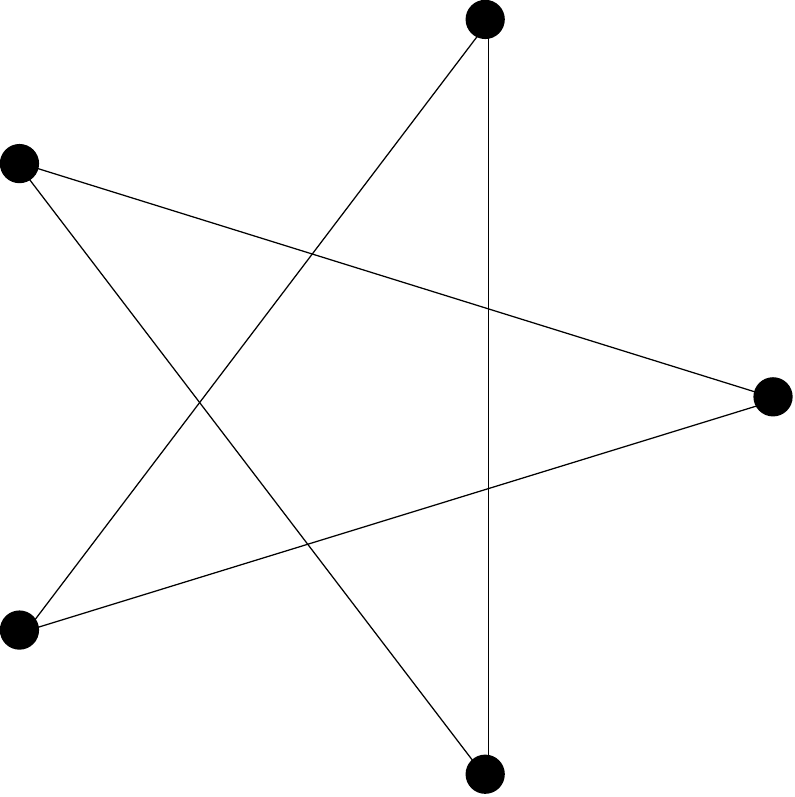}
\caption{Regular pentagon and pentagram.}
\end{figure}

\medskip{}




\section{Buffon transformation for polyhedra}

Recall first some basic notions of graph theory and the relation with polyhedra.

A {\it graph} $\Gamma=(\mathcal V,\mathcal E)$ consists of a finite set $\mathcal V$ (vertices), together with
a subset $\mathcal E\subseteq \mathcal V\times \mathcal V$ (edges). We will assume that the graph has no loops $[i,i], \,\, i \in \mathcal V$
and is {\it undirected} which means that
for each edge $[i,j]\in \mathcal E$ we also have $[j,i]\in \mathcal E$.

We say that the vertices $i$ and $j$ are {\it adjacent} 
and write $i \sim j$ if there is an edge $[i,j] \in \mathcal E$ connecting them.
The {\it degree} $d_i$ of a vertex $i$ is the number of the adjacent vertices. 

A graph is {\it connected} when there is a path between any two vertices.
A graph is called {\it 3-connected} if for every pair of vertices $i$ and $j$ there
are at least three paths from $i$ to $j$, whose only vertices (or
edges) in common are $i$ and $j$.

A graph is called {\it planar} if an isomorphic copy of the
graph can be drawn in a plane, so that the edges which
join the vertices only meet (intersect) at vertices.  

For every polyhedron $P$ one can consider the {\it 1-skeleton} $\Gamma (P)$, which is
the graph formed by the vertices and edges of $P.$

One of the oldest results in polytope theory is a remarkable theorem
by Ernst Steinitz. It is often referred to as the Steinitz'
{\it fundamental theorem of convex polyhedra} and gives a completely
combinatorial characterization of the graphs $\Gamma,$  
which can be realised as 1-skeletons of 3-dimensional
polytopes (see \cite{Grun,Z}). 

\begin{thm}
{\bf (Steinitz, 1922)} A graph $\Gamma$ is isomorphic to the 1-skeleton
of a 3-dimensional convex polyhedron P if and only if $\Gamma$ is
planar and 3-connected.\medskip{}
\end{thm}

The proof given by Steinitz uses a combinatorial reduction technique.
A sequence of transformations of $\Gamma$ into simpler graphs lead
to the tetrahedral graph $K_{4}$. Reversing the order of these operations
one obtains a polyhedral realization of the original graph $\Gamma$.

A graph is called {\it regular} when every graph vertex has the same
degree.



Let $P$ be a simplicial polyhedron in $\mathbb R^3$ with vertices $r_1, \dots, r_n.$
Define the {\it Buffon transformation} $B(P)$ as a new polyhedron with the vertices being the centroids 
of all edges, which meet at a vertex \cite{VW,W2}:
\begin{equation}
\label{BT}
B(r_i)=\sum_{j \sim i} \frac{1}{2d_i}(r_i+r_j),
\end{equation}
where $d_i$ is the degree of the vertex $r_i.$

Consider also the linear {\it Buffon operator} $B:\mathcal F(\mathcal V) \rightarrow \mathcal F(\mathcal V),$ where $\mathcal F(\mathcal V)$ is the vector space of functions on the vertices 
of the graph $\Gamma = \Gamma(P),$ defined by the same formula:
\begin{equation}
\label{BTop}
B(f)(i)=\sum_{j \sim i} \frac{1}{2d_i}(f(i)+f(j)),\,\, f \in \mathcal F(\mathcal V).
\end{equation}




{\it Remark.} One can define the Buffon transformation $B_F$ by taking  the centroids of
the centroids of all the faces meeting at a vertex \cite{VW,W2}, 
but for simplicial polyhedra $P$ we have a simple relation for the corresponding operators
$$
B_F=\frac{4}{3}B - \frac{1}{3} I,
$$
which means that the result of the Buffon procedure on faces will be the same as the one on edges.

The matrix of the Buffon transformation in a natural basis in $\mathcal F(\mathcal V)$ has the form
\begin{equation}
\label{BTmat}
B=\frac{1}{2}(I+D^{-1}A)=\frac{1}{2}(I+P),
\end{equation}
where $A$ is the {\it adjacency matrix}: 
$A_{ij}=1$ if $i$ is adjacent to $j$ and 0 otherwise,
$D$ the diagonal matrix with the degrees of vertices
$d_{i}$ on the diagonal, and $P$ is the matrix of {\it transition probabilities} of the Markov chain describing the  random walk on graph $\Gamma:$
$P_{ij}= 1/d_j$ when $j$ is adjacent to $i$ and 0 otherwise (see \cite{L3}).

Note that unless $\Gamma$ is a regular
graph, matrix $B$ is not symmetric. 
In order to bring it to a symmetric form we introduce the {\it normalized
adjacency matrix}
\begin{equation}
\label{normad}
N=D^{-\frac{1}{2}}AD^{-\frac{1}{2}}
\end{equation}
with matrix elements $N_{ij}= 1/\sqrt{d_id_j}$ if $i$ is adjacent to $j$ and 0 otherwise.
It is easy to see that
\[
B=\frac{1}{2}(I+D^{-\frac{1}{2}}ND^{\frac{1}{2}})=\frac{1}{2}D^{-\frac{1}{2}}(I+N)D^{\frac{1}{2}},
\]
so $B$ is conjugated to the symmetric matrix
$\tilde B=1/2(I+N).$

In particular, this means that all the eigenvalues of $B$ are real. The maximal eigenvalue is $\lambda_0=1$ and the corresponding eigenvector is $(1,\dots,1)^T.$

Now we ask the same question: what is the limiting shape of $B^n(P)$ when $n$ goes to infinity ?

By the same arguments using the Subspace Iteration Theorem the answer is given by the subdominant 
eigenspace of the corresponding Buffon operator $B.$ In general it is one-dimensional, which means that the limiting shape is one-dimensional.
However, if we assume additional symmetry we have a three-dimensional limiting shape. 
To see this we need some results from spectral graph theory, which we present in the next section.

\section{Colin de Verdi\`ere invariant and null space representation}

In 1990 Yves Colin de Verdi\`ere \cite{CdV} introduced a new spectral
graph invariant $\mu(\Gamma).$ Roughly speaking, $\mu(\Gamma)$ 
is the maximal multiplicity of the second largest eigenvalue of the matrices $C$ with the property
$C_{ij}=C_{ji}>0$ for adjacent $i$ and $j,$ $C_{ij}=0$ for non-adjacent $i$ and $j$ and arbitrary diagonal elements $C_{ii}.$
The precise definition is as follows.

Let $\Gamma$ be a connected undirected graph with the vertex set
$\left\{ 1,\ldots,n\right\} $. Let $\mathcal{M}_{\Gamma}$ denote
the set of symmetric matrices $M=(M_{ij})\in\mathbb{R}^{V\times V}$
associated with $\Gamma$ satisfying
\begin{enumerate}
\item $M_{ij}\begin{cases}
<0, & ij\in E\\
=0, & ij\notin E
\end{cases}$;
\item $M$ has exactly one (simple) negative eigenvalue.
\end{enumerate}
\medskip{}

$M$ is said to satisfy the {\it Strong Arnold Property} if the relation $MX=0$ with a symmetric
$n\times n$ matrix $X$ such that $X_{ij}=0$ for any adjacent
$i$ and $j$ and for $i=j$ implies that $X=0$. This property
is a restriction, which excludes some degenerate choices of the edge
weights and the diagonal entries.

The {\it Colin de Verdi\`ere invariant} $\mu\left(\Gamma\right)$ is
the largest corank of matrices from the set $\mathcal{M}_{\Gamma}$
satisfying the Strong Arnold Property. A matrix $M\in\mathcal{M}_{\Gamma}$
with corank $\mu\left(\Gamma\right)$ is called a {\it Colin de Verdi\`ere
matrix} of $\Gamma$.

After the change of sign and shift by a scalar matrix $C=cI-M$ the corank, which is the dimension of the null space of
$M$ becomes the multiplicity of the second largest eigenvalue of $C.$. 

Colin de Verdi\`ere characterised all the graphs with parameter $\mu\left(\Gamma\right)\leq 3.$ 

A graph is called {\it outerplanar} if it can be drawn in the plane without
crossings in such a way that all of the vertices belong to the unbounded
face of the drawing. 

\begin{thm} {\bf (Colin de Verdi\`ere, 1990)}
\begin{itemize}
\item $\mu\left(\Gamma\right)\leq1$ if and only if $\Gamma$ is a path;
\item $\mu\left(\Gamma\right)\leq2$ if and only if $\Gamma$ is outerplanar;
\item $\mu\left(\Gamma\right)\leq3$ if and only if $\Gamma$ is planar.\label{nu-planarity}
\end{itemize}
\end{thm}

The planarity characterization is a remarkable result, which will be important for us. 
The "only if" part is relatively simple and follows from Kuratowski's characterisation of the planar graphs \cite{Harari}.
The original proof of the "if" part was quite involved.
Van der Holst \cite{Holst1} substantially simplified it and showed that for 3-connected planar graphs
the Strong Arnold property does not play any role. 

\begin{corollary} {\bf (Van der Holst, 1995)}\label{corH}
For any matrix $M$ from $\mathcal M_{\Gamma}$ the corank of $M$ can not be larger than 3.
\end{corollary}

In \cite{L1} Lov\'asz found an explicit way of constructing the Colin de Verdi\`ere matrix 
for any 3-connected planar graph $\Gamma$ using the Steinitz' realisation of $\Gamma$ as 
a 1-skeleton of a convex polyhedron $P.$ 
This result will be crucial for us, so we will sketch here the main steps of his construction following \cite{L1}.

Recall first the notion of polarity for polyhedra in $\mathbb R^3$, see e.g. \cite{Z}.
Let $P$ be any convex polytope in $\mathbb{R}^{3}$, containing the origin
in its interior. The {\it polar polyhedron} $P^*$ is defined as
$$
P^*=\{ y \in \mathbb R^3: (y,x) \leq 1 \,\, {\text for \,all} \,\, x \in P\},
$$
where $( , )$ denote the scalar product in $\mathbb R^3.$
It is known that $P^*$ is also a convex polyhedron and 
the 1-skeleton of $P^*$ is the planar dual graph $\Gamma^*=(\mathcal V^*, \mathcal E^*)$ with vertices 
corresponding to the faces of $P$ and edges corresponding to edges of $P$ \cite{Z}.

Now let $P \subset \mathbb R^3$ be Steinitz' realisation of graph $\Gamma,$ so that $\Gamma$ is isomorphic to 1-skeleton $\Gamma(P).$
We can always assume that $P$ contains the origin inside it.
Consider its polar polyhedron $P^*$

Let $u_{i}$ and $u_{j}$ be two adjacent vertices of $P$, and 
$w_{f}$ and $w_{g}$ be the endpoints of the corresponding edge of
$P^{*}$. Then by the definition of polarity we have
\[
(w_{f},u_{i})=(w_{g},u_{i})=1.
\]
 This implies that $w_{f}-w_{g}$ is perpendicular to $u_i$, and similarly to $u_j.$
Hence the vectors $w_{f}-w_{g}$ and the cross-product $u_{i}\times u_{j}$ are parallel
and we can find the coefficients $M_{ij}$ such that
\[
w_{f}-w_{g}=M_{ij}(u_{i}\times u_{j}).
\]
We can always choose the labelling of $w_{f}$ and $w_{g}$ in such a way that $M_{ij}<0$.

This defines $M_{ij}$ for adjacent $i\neq j.$ For non-adjacent $i$ and $j$ we define $M_{ij}$ to be zero.
To define $M_{ii}$
consider the vector
\[
u_{i}^{'}=\sum_{j \sim i}M_{ij}u_{j}.
\]
Then
\[
u_{i}\times u_{i}^{'}=\sum_{j \sim i}M_{ij}u_{i}\times u_{j}=\sum(w_{f}-w_{g}),
\]
where the last sum is taken over all edges $fg$ of the face
of $P^{*}$ corresponding to $i$, oriented counterclockwise. Since
this sum is zero we have
\[
u_{i}\times u_{i}^{'}=0,
\]
which means that $u_{i}$ and $u_{i}^{'}$ are parallel. Therefore we can define $M_{ii}$ by the relation
\[
u_{i}^{'}=-M_{ii}u_{i}.
\]

\begin{thm} {\bf (Lov\'asz, 2000)}
The matrix $M$ described above is a Colin de Verdi\`ere matrix
for the graph $\Gamma$.
\end{thm}

Indeed by construction $M$ has the right pattern of zeros
and negative elements. The condition $u_{i}^{'}=-M_{ii}u_{i}$ can
be written in the following form 
\[
\sum_{j}M_{ij}u_{j}=0.
\]
This means that each coordinate of the $u_{i}$ defines a vector in
the kernel of $M$ and hence $M$ has corank at least 3.
But by Corollary \ref{corH} it can not be larger than 3, so the corank is 3 and thus maximal.

To prove that $M$ has exactly one negative eigenvalue one can use the classical Perron-Frobenius theorem, see e.g. \cite{G}.

\begin{thm}
{\bf (Perron-Frobenius, 1912)} If a real matrix has non-negative entries
then it has a nonnegative real eigenvalue $\lambda$ which has maximum
absolute value among all eigenvalues. This eigenvalue $\lambda$ has
a real eigenvector with non-negative coordinates. If the matrix is irreducible, then $\lambda$ has multiplicity
1 and the corresponding eigenvector can be chosen to be positive.
\end{thm}

Choosing sufficiently large $c>0$ we have the matrix $cI-M,$ which has 
non-negative entries and irreducible, so we can apply the Perron-Frobenius
Theorem to conclude that the smallest eigenvalue of $M$ has multiplicity 1. It must be negative since we know that
 the eigenvalue $0$ has multiplicity at least 3. The fact that there are no more negative multiplicities requires a considerable work using 
the connectivity of the space of Steinitz' realisations, see \cite{L1}.

Conversely, having a Colin de Verdi\`ere matrix $M \in \mathcal M_{\Gamma}$ one can consider the following {\it null space representation} 
$\nu:  \mathcal V =\{1,2, \dots, n\} \rightarrow \mathbb R^3$ (see \cite{L2}).
 
Choose a basis $a_1, a_2, a_3$ in the kernel of $M$ and consider a $3 \times n$ matrix $X$ with rows being the coordinates of $a_1, a_2, a_3.$
Then the columns $u_i, \, i=1, \dots, n$ of this matrix give the set of 3-vectors, defining the map $\nu.$ 
The problem is that in general they will not be vertices of a convex polyhedron, but Lovasz \cite{L1} showed that after some scaling 
$u_i \rightarrow \mu_i u_i$ this is the case (such a scaling he called {\it proper}). At the level of the Colin de Verdi\`ere matrices this corresponds to the change $M \rightarrow DMD$, where $D=diag\, (\mu_1, \dots, \mu_n)$ is a non-degenerate diagonal matrix, which obviously preserves the properties of $\mathcal M_{\Gamma}.$ 

\begin{thm} {\bf (Lov\'asz, 2000)}
For a 3-connected planar graph $\Gamma$ any Colin de Verdi\`ere matrix $M \in \mathcal M_{\Gamma}$ can be properly scaled, 
so that null space representation gives a convex polyhedron with 1-skeleton isomorphic to $\Gamma$.
\end{thm}

Note that the change of basis in the kernel of $M$ corresponds to a linear transformation of $\mathbb R^3$, so the corresponding 
polyhedron is defined only modulo affine transformation.

Now we are ready to prove our main result.

\section{Proof of the main theorem}

Let $G$ be a Platonic group and $\Gamma$ a $G$-invariant planar
3-connected graph. 

We know after Steinitz that $\Gamma$ can be realized by a 3-dimensional
convex polyhedron $P$, but in the presence of symmetry Mani \cite{Mani} 
showed that there is a symmetric realisation $P_G \subset \mathbb R^3.$

\begin{thm} {\bf (Mani, 1971)}
There exists a convex polyhedron $P_G \subset \mathbb R^3$ 
with the group of isometries isomorphic to $G$ and with 1-skeleton isomorphic to $\Gamma.$  
\end{thm}

Since $\Gamma$ is
planar and 3-connected, its Colin de Verdi\`ere invariant $\mu(\Gamma)$
must be 3. Let $M$ be the Colin de Verdi\`ere matrix given by Lov\'asz construction applied 
to Mani's version of Steinitz realisation $P_G.$

Let $N$ be the normalised adjacency matrix (\ref{normad}). We know that 
the matrix of the Buffon transformation $B$ is related to $N$ by
\[
B=\frac{1}{2}D^{-\frac{1}{2}}(I+N)D^{\frac{1}{2}}
\]
and that its largest eigenvalue is $\lambda_0=1.$ Let $\lambda_1$ be the second largest eigenvalue of $B.$ 
We would like to show that it has multiplicity 3.

To do this consider the symmetric matrix
\[
\widehat{B}=-\frac{1}{2}N+(\lambda_{1}-\frac{1}{2})I.
\]
It is easy to see that
$\widehat{B}\in\mathcal{M}_{\Gamma}$ and that the corank of $\widehat{B}$
is precisely the multiplicity of $\lambda_1.$ 

Define a parameter family of matrices 
\begin{equation}
M_{t}=(1-t)M-t\widehat{B}, \quad t\in[0,1]
\end{equation}
where $M$ is the Colin de Verdi\`ere matrix defined above.

Since $M_{t}$ is $G$-invariant,
the group $G$ acts on the kernel of $M_{t}$. 
When $t=0$ we know that the kernel of $M(0)=M$ has dimension 3 and by Lov\'asz result \cite{L1}
the corresponding representation of $G$ is standard geometric 
by the isometries of $P_G.$

Since this representation is irreducible and the set of 3-dimensional representations of $G$ is discrete, 
by continuity arguments the kernel will remain 3-dimensional geometric representation
for all $t\in [0,1]$, in particular for $t=1$. 

These arguments will not work only if $0$ collides with another eigenvalue. But this could not happen 
with the negative eigenvalue because of the Perron-Frobenius theorem. 
In particular, all matrices $M_t$ belong to $\mathcal{M}_{\Gamma}.$ 
If this happens with a positive eigenvalue we will have the corank of the corresponding $M_t$ to be at least 4, 
which contradicts to the Colin de Verdi\`ere result. 

Thus we have proved that the kernel of $M_1=-\widehat{B}$ is 3-dimensional, 
and hence the same is true for the subdominant eigenspace of the Buffon operator $B.$
The limiting shape is given essentially by the null space representation construction, 
but the proper scaling may not hold. However, the very existence of a proper scaling \cite{L1, L2}
and the assumption of simplicity imply that the corresponding vectors $u_i$ are the vertices of a 
certain star-shaped polyhedron with 1-skeleton isomorphic to $\Gamma.$
The triakis tetrahedron example below shows that the proper scaling is indeed not automatic, 
so the convexity property does not necessary hold. 
 
This completes the proof of Theorem \ref{thm:main}.

\section{Concluding remarks.}

The Buffon regularisation procedure can be interpreted as search of an ideal shape of a given polyhedron 
and in that sense can be considered as one of the earliest examples of the trend, 
popular in modern differential geometry. 

For manifolds this usually leads to the solutions of certain nonlinear PDEs like the
mean curvature flow in the theory of minimal surfaces \cite{Huisken} or the celebrated Ricci flow in Thurston's geometrization programme \cite{MT}. 
Our case is conceptually closer to the description of the minimal submanifolds in the unit sphere 
using the eigenfunctions of the Laplace-Beltrami operator, see \cite{KN, T}.

The main difference with the differential case is that the generic graphs are 
much less regular objects than manifolds, even under our assumption of Platonic symmetry.
The crucial thing here is a large multiplicity of the second eigenvalue of the Buffon operator. 
How to guarantee this is a good question.

The symmetry assumption seems to be natural. In this relation we would like to mention 
an interesting result of  Mowshowitz \cite{Mowsh}, 
who showed that if all eigenvalues of the adjacency
matrix $A$ of a graph are different, then every automorphism of $A$
has order 1 or 2. Some interesting related results for the graphs with vertex transitive group action can be found in \cite{IP}. Note that in our case the group action is far from being vertex transitive.

An interesting question concerns the decomposition of $\mathcal F(\mathcal V)$ into the irreducible $G$-modules with respect to the Buffon spectrum. We saw that the geometric representation always appears at the subdominant level, but we do not know much about higher level. For the regular polyhedra the answer is given in Appendix A.

It would be interesting to understand what our geometric analysis means for related random walk on the corresponding graphs.

Finally, a natural question is what happens in higher dimension. We believe that for the simplicial polyhedra 
 we should expect similar result if we assume the symmetry under an irreducible  Coxeter group. 
 Note that in dimension 4 we have 6 regular polyhedra with the symmetry groups $A_4=S_5$, $B_4$, $F_4$ and $H_4,$
 while in dimension more than 4 we have only analogues of tetrahedron, cube and octahedron.

\section{Acknowledgements.}

We are grateful to Jenya Ferapontov, Steven Kenny, Boris Khesin and L\'aszl\'o Lov\'asz for very helpful discussions.
Special thanks are to Graham Kemp, who was part of these discussions for quite a while.

 \bigskip
 
{\bf Appendix A. The symmetry groups of Platonic solids and their characters}

\medskip

The symmetry group of a regular tetrahedron is $S_4$ and is isomorphic to the permutation group of the vertices.

The full symmetry group of the octahedron is the same as for the cube: $G=S_4 \times \mathbb Z_2$.
$S_4$ is the rotation subgroup, which is isomorphic to the permutation group of the 4 long diagonals, and 
$\mathbb Z_2$ corresponds to the central symmetry of cube. 

For the icosahedron and dodecahedron the full symmetry group  is known to be $A_5 \times \mathbb Z_2,$ 
where $A_5 \subset S_5$ is the alternating subgroup of $S_5$ describing the rotational symmetry and 
$\mathbb Z_2$ is again the central symmetry of the solids.

The irreducible representations of the group $G=H \times \mathbb Z_2$ have the the form $V_1\otimes V_2$, where 
$V_1$ and $V_2$ are irreducible representations of $H$ and $\mathbb Z_2$ respectively. Note that $V_2$ is either trivial or sign representation of $\mathbb Z_2,$ which we will denote respectively by 1 and $\varepsilon$.
Thus we need only the character tables of the groups $S_4$ and $A_5$, which in the notations of Fulton and Harris \cite{FH} are given below in Tables \ref{Table:charssss} and \ref{Table:charaaaaa}.

\begin{table}[h]
	\centering
		\begin{tabular}{|r|ccccc|}
		\hline
		  24 & 1 & 6& 8 & 6 &3 \\
			$S_4$&  1 &(12)&(123) & (1234) & (12)(34) \\
			\hline 
			$U$ & 1 &1 & 1 &1 & 1\\
			$U'$ & 1 & -1 & 1 &-1 & 1\\
			$V$ & 3 & 1 & 0 &-1 & -1\\
			$V'$ & 3 & -1 & 0 & 1 & -1 \\
			$W$ & 2 & 0 &-1 & 0 &2 \\
			\hline
		\end{tabular}
				\caption{The character table of $S_4$.}\label{Table:charssss}
\end{table}

\begin{table}[h]
	\centering
		\begin{tabular}{|r|ccccc|}
		\hline
		  60 & 1 & 20&  15 & 12 											&12 \\
			$A_5$&  1 &(123)&(12)(34) & (12345)				 & (21345) \\
			\hline 
			$U$ & 1 &1  & 1 &1												 & 1\\
			$V$ & 4 &1  & 0 &-1 											 & -1\\
			$W$ & 5 & -1& 1 &0 												 & 0\\
			$Y$ & 3 & 0 & -1&  $\frac{1+ \sqrt{5}}{2}$ & $\frac{1- \sqrt{5}}{2}$ \\
			$Z$ & 3 & 0 &-1 & $\frac{1- \sqrt{5}}{2}$  & $\frac{1+ \sqrt{5}}{2}$ \\
			\hline
		\end{tabular}
		\caption{The character table of $A_5$.}\label{Table:charaaaaa}
\end{table}

With these notations the geometric representations are: $V$ for tetrahedral group $G=S_4,$ 
$\varepsilon V'=V'\otimes \varepsilon$ for cube/octahedral group $G=S_4\times \mathbb Z_2$ and $\varepsilon Y=Y\otimes \varepsilon$ for icosahedral/dodecahedral group $G=A_5\times \mathbb Z_2.$

The corresponding decompositions of the space of functions on the vertices into irreducible $G$-modules are
\begin{equation}
\label{T}
\mathcal F(T)=U\oplus V
\end{equation}
for tetrahedron,
\begin{equation}
\label{O}
\mathcal F(O)= U \oplus \varepsilon V' \oplus W
\end{equation}
 for octahedron,
\begin{equation}
\label{C}
\mathcal F(C)= U \oplus \varepsilon V'  \oplus V  \oplus \varepsilon U'
\end{equation}
 for cube,
\begin{equation}
\label{I}
\mathcal F(I)= U \oplus \varepsilon Y \oplus W \oplus \varepsilon Z
\end{equation}
 for icosahedron,
\begin{equation}
\label{D}
\mathcal F(D)= U \oplus \varepsilon Y  \oplus W \oplus \varepsilon V \oplus  V \oplus \varepsilon Z
\end{equation} for dodecahedron.

We have ordered them according to the appearance in the spectrum of the Buffon operator. It turns out that in all these cases
the spectral decomposition coincides with $G$-decomposition (see the examples below). Note that the first two are always trivial and geometric representations in agreement with our result.

\bigskip
\bigskip
\bigskip

{\bf Appendix B. Examples of Buffon realizations of polyhedra}

\bigskip

For the polyhedra $P$ with combinatorial structure of Platonic solids 
the Buffon procedure leads to the polyhedron $P_B$, 
which is affine equivalent to the regular realisation of $P.$

Since in the regular case the Buffon matrix $B$ can be replaced by the adjacency matrix $A$ 
the calculations are essentially the same as in \cite{MC}, 
where one can find a lot more experimental data.
The calculation of spectra of regular polytopes can be found in \cite{ST}.

We present here the most instructive examples of Buffon realisations of regular, Archimedean and Catalan solids.
All the calculations and pictures were made using Mathematica.

Recall that the Archimedean solids (also referred to as the semi-regular polyhedra) are the convex polyhedra with faces being regular polygons of two or more different types arranged in the same way about each vertex. Solids with a dihedral group of symmetries (e.g., regular prisms and
antiprisms) are not considered to be Archimedean solids. With this
restriction there are 13 Archimedean solids. 
For Archimedean solids the affine $B$-regular version is in general is not affine equivalent to the standard one (see below the example of truncated cube).

The Catalan solids are duals of the Archimedean solids. The Catalan solids
are convex polyhedra with regular vertex figures (of different types) and with equal dihedral angles.
For Catalan solids the affine $B$-regular versions may not be convex
or, in non-simplicial case, may even not exist (see the examples below).

We start with the regular cases of icosahedron and dodecahedron 
to show the relation with $G$-decomposition and to look at the embeddings 
related to other eigenvalues.

\bigskip{}


{\bf The Icosahedron}

\begin{figure} [H]
   \centering
  \includegraphics[scale=0.8]{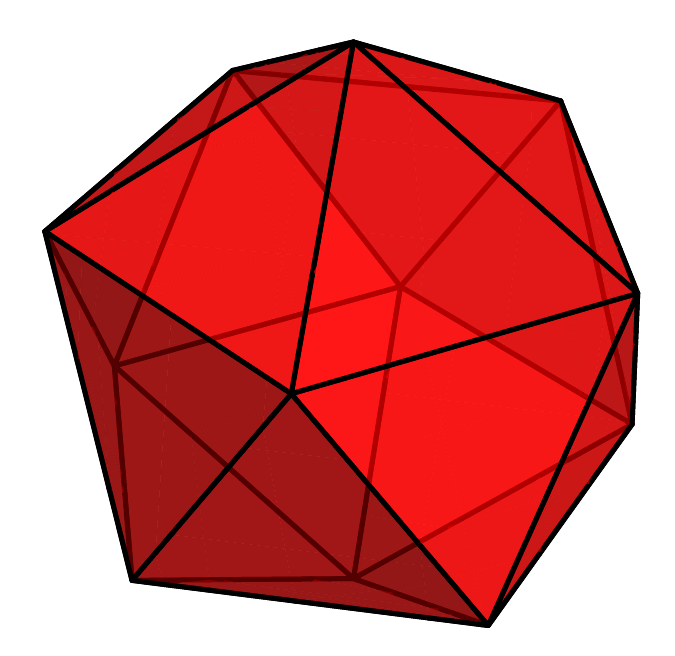}\quad{}\includegraphics[scale=0.8]{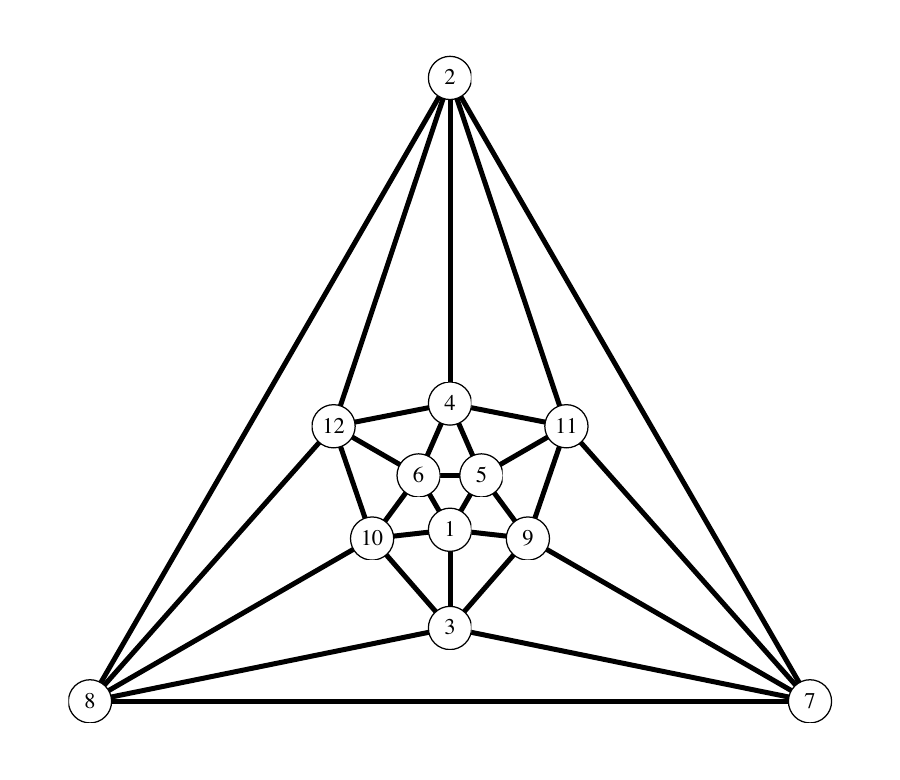}
  \medskip
   \caption{The icosahedron and its $1$-skeleton graph.}
\end{figure}

The corresponding Buffon spectrum is:
\[
\left\{ 1^{(1)},\frac{1}{10}\left(5+\sqrt{5}\right)^{(3)},\frac{2}{5}^{(5)},\frac{1}{10}\left(5-\sqrt{5}\right)^{(3)}\right\} 
\]
in agreement with (\ref{I}).

The eigenspaces corresponding to the second highest eigenvalue
and its conjugate eigenvalue $\frac{1}{10}\left(5-\sqrt{5}\right)$
are:
\[
X_{2}=\left(\begin{array}{c}
a(\alpha+\gamma)-\beta\\
\beta-a(\alpha+\gamma)\\
a(\beta-\gamma)-\alpha\\
a(\gamma-\beta)+\alpha\\
-\gamma\\
a(\beta-\alpha)-\gamma\\
a(\alpha-\beta)+\gamma\\
\gamma\\
-\alpha\\
-\beta\\
\beta\\
\alpha
\end{array}\right)\quad X_{4}=\left(\begin{array}{c}
b(\alpha+\gamma)-\beta\\
\beta-b(\alpha+\gamma)\\
b(\beta-\gamma)-\alpha\\
\alpha+b(\gamma-\beta)\\
-\gamma\\
b(\beta-\alpha)-\gamma\\
b(\alpha-\beta)+\gamma\\
\gamma\\
-\alpha\\
-\beta\\
\beta\\
\alpha
\end{array}\right),
\]
where $a=\frac{1}{2}(1-\sqrt{5})$ and $b=\frac{1}{2}(1+\sqrt{5})$.

Geometrically, $X_{2}$ describes an affine regular icosahedron, while $X_{4}$
corresponds to an affine great icosahedron, which is one of four Kepler-Poinsot regular star polyhedra.

\begin{figure} [H]
   \centering
  \includegraphics[scale=0.7]{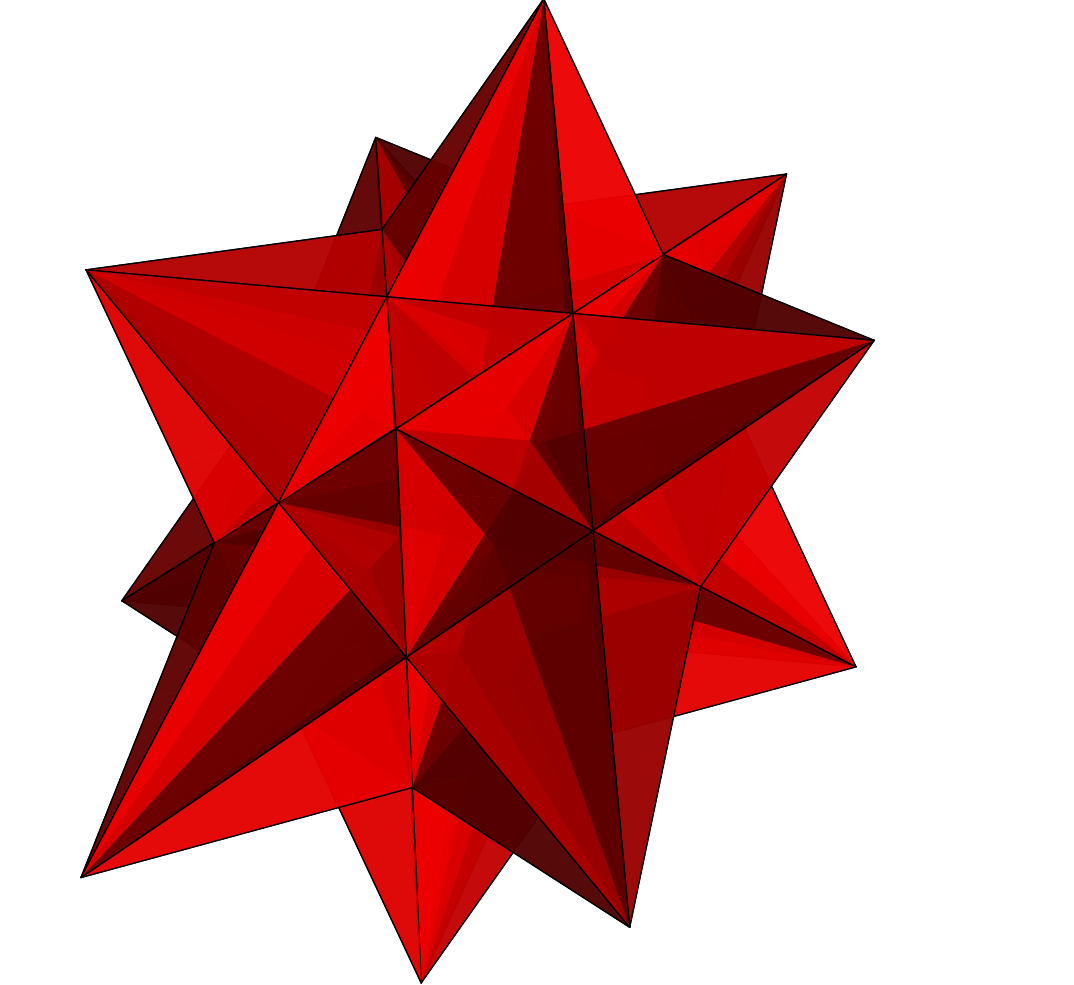}
     \medskip
   \caption{The affine great icosahedron.}
 \end{figure}
 
The eigenspace $X_3$ describes the 5-dimensional realisation of an icosahedron as a 5-simplex:
6 pairs of opposite vertices identified with 6 vertices of the simplex.

\bigskip


{\bf The Dodecahedron}

\medskip

 \begin{figure} [H]
   \centering
  \includegraphics[scale=0.55]{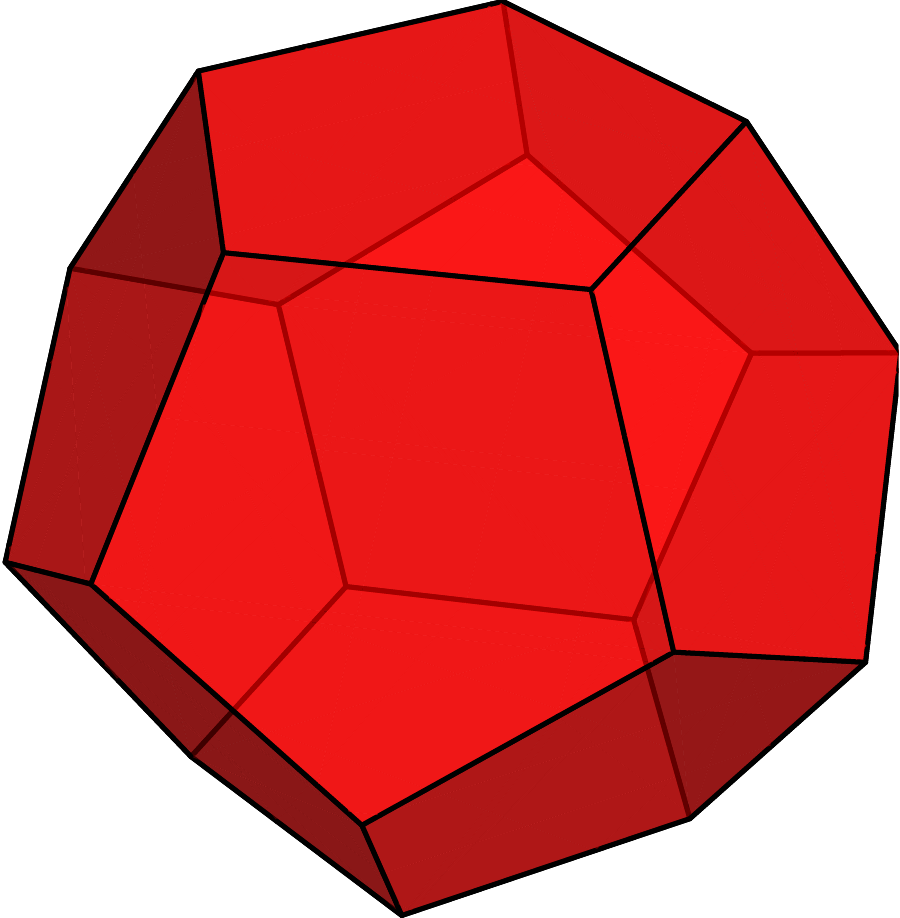}\quad{}\includegraphics[scale=0.8]{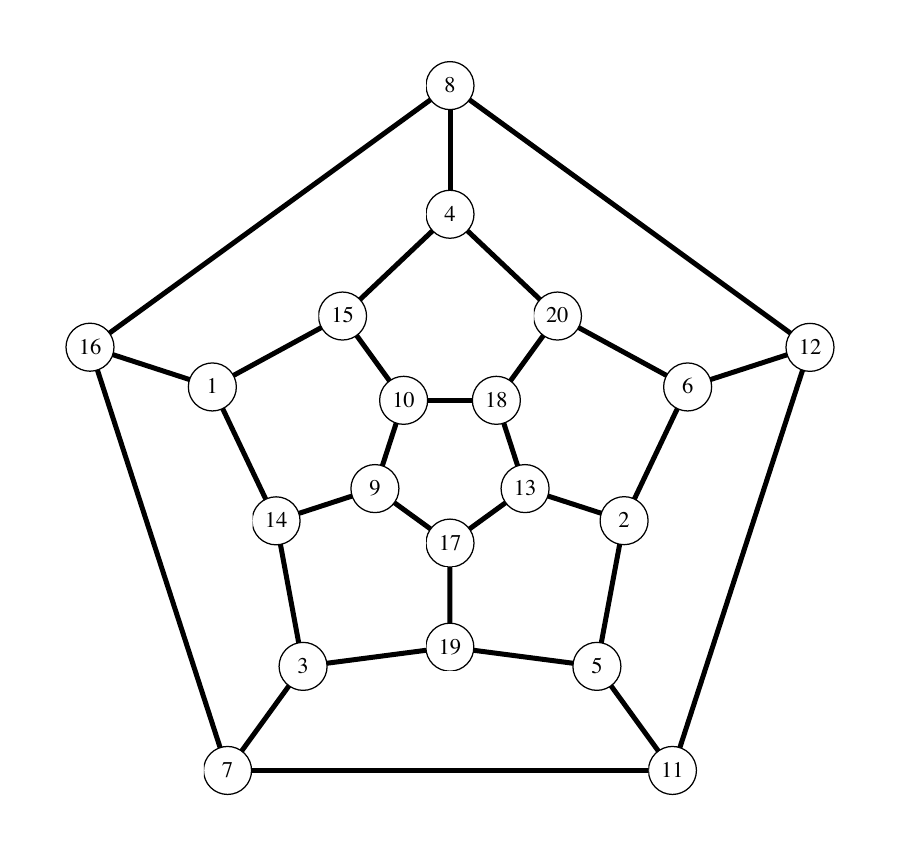}
    \medskip
   \caption{The dodecahedron and its $1$-skeleton graph.}
\end{figure}

The corresponding Buffon spectrum is 
\[
\left\{ 1^{(1)},\frac{1}{6}\left(3+\sqrt{5}\right)^{(3)},\frac{2}{3}^{(5)},\frac{1}{2}^{(4)},\frac{1}{6}^{(4)},\frac{1}{6}\left(3-\sqrt{5}\right)^{(3)}\right\} 
\]
in agreement with (\ref{D}).

The eigenspaces corresponding to the second highest eigenvalue $\lambda_{2}=\frac{1}{6}\left(3+\sqrt{5}\right)$ and to its conjugate $\lambda_{6}=\frac{1}{6}\left(3-\sqrt{5}\right)$  are:

\[
X_{2}=\left(\begin{array}{c}
-a\gamma-2\alpha-b\beta\\
a\gamma+2\alpha+b\beta\\
-\alpha\\
-\beta\\
b(\alpha+\beta)-\gamma\\
\sqrt{5}\alpha+\beta-\gamma\\
-\gamma\\
a(\beta-\alpha)-\gamma\\
a\beta-2\alpha+b\gamma\\
b(\gamma-\alpha)-\beta\\
b(\alpha-\gamma)+\beta\\
-a\beta+2\alpha-b\gamma\\
\beta-a(\alpha+\gamma)\\
-\sqrt{5}\alpha-\beta+\gamma\\
\gamma-b(\alpha+\beta)\\
a(\alpha+\gamma)-\beta\\
a(\alpha-\beta)+\gamma\\
\gamma\\
\beta\\
\alpha
\end{array}\right)\quad X_{6}=\left(\begin{array}{c}
-a\beta-2\alpha-b\gamma\\
a\beta+2\alpha+b\gamma\\
-\alpha\\
-\beta\\
a(\alpha+\beta)-\gamma\\
-\sqrt{5}\alpha+\beta-\gamma\\
-\gamma\\
b(\beta-\alpha)-\gamma\\
a\gamma-2\alpha+b\beta\\
a(\gamma-\alpha)-\beta\\
a(\alpha-\gamma)+\beta\\
-a\gamma+2\alpha-b\beta\\
\beta-b(\alpha+\gamma)\\
\sqrt{5}\alpha-\beta+\gamma\\
\gamma-a(\alpha+\beta)\\
b(\alpha+\gamma)-\beta\\
b(\alpha-\beta)+\gamma\\
\gamma\\
\beta\\
\alpha
\end{array}\right),
\]
where $a=\frac{1}{2}(1-\sqrt{5})$ and $b=\frac{1}{2}(1+\sqrt{5})$.

Geometrically, $X_{2}$ describes an affine regular dodecahedron, while
$X_{6}$ corresponds to an affine version of the great stellated dodecahedron, which is another Kepler-Poinsot polyhedron (see Fig. 7).

It is a bit puzzling that the remaining two Kepler-Poinsot polyhedra (small stellated dodecahedron and great dodecahedron) 
seem to not appear in Buffon approach.

\begin{figure} [H]
   \centering
  \includegraphics[scale=0.7]{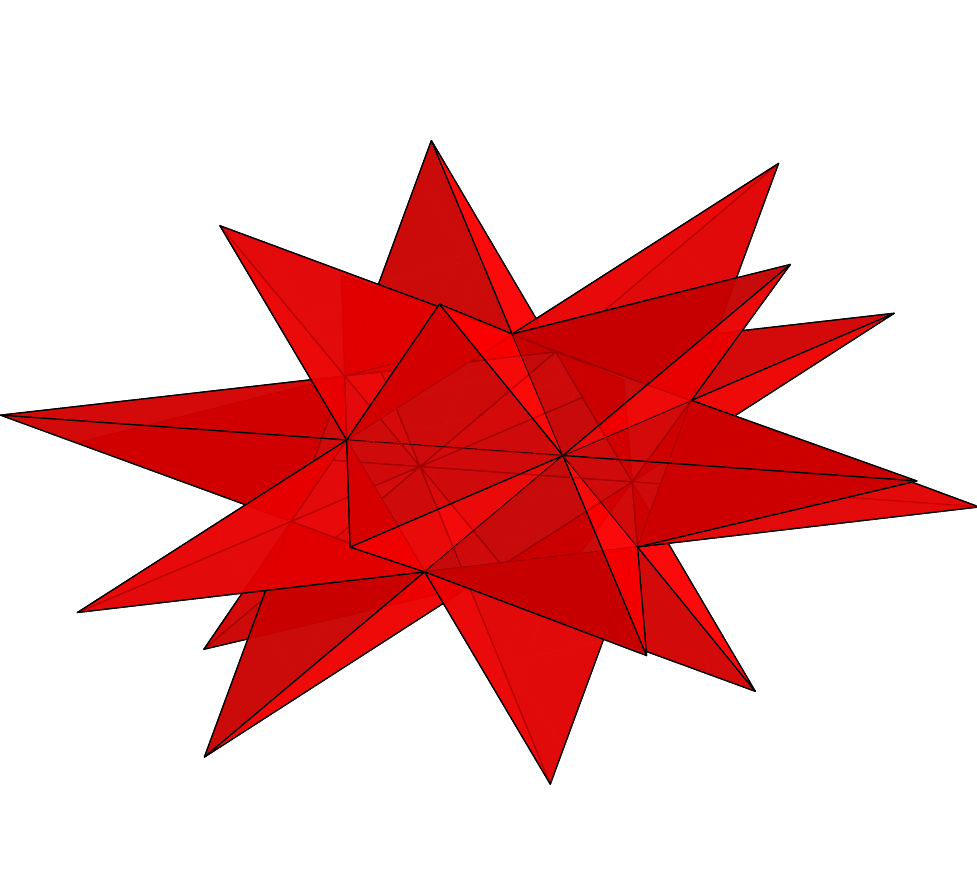}
  \medskip
  \caption{Affine great stellated dodecahedron.}
\end{figure}

The 3 remaining eigenspaces are: 

\[X_3=\left(
\begin{array}{c}
 -\gamma -\delta +\rho  \\
 -\gamma -\delta +\rho  \\
 \alpha  \\
 \beta  \\
 -\alpha +\beta -\delta  \\
 \alpha -\beta -\gamma  \\
 \gamma  \\
 \delta  \\
 -\beta +\delta -\rho  \\
 -\alpha +\gamma -\rho  \\
 -\alpha +\gamma -\rho  \\
 -\beta +\delta -\rho  \\
 \rho  \\
 \alpha -\beta -\gamma  \\
 -\alpha +\beta -\delta  \\
 \rho  \\
 \delta  \\
 \gamma  \\
 \beta  \\
 \alpha 
\end{array}
\right), 
\quad
X_{4}=\left(
\begin{array}{c}
 \gamma +\delta  \\
 -\gamma -\delta  \\
 -\alpha  \\
 -\beta  \\
 \alpha -\delta  \\
 \beta -\gamma  \\
 -\gamma  \\
 -\delta  \\
 \alpha -\gamma -\delta  \\
 \beta -\gamma -\delta  \\
 -\beta +\gamma +\delta  \\
 -\alpha +\gamma +\delta  \\
 -\alpha -\beta +\gamma +\delta  \\
 \gamma -\beta  \\
 \delta -\alpha  \\
 \alpha +\beta -\gamma -\delta  \\
 \delta  \\
 \gamma  \\
 \beta  \\
 \alpha 
\end{array}
\right),
\quad 
X_{5}=\left(
\begin{array}{c}
 2 \alpha +2 \beta +\gamma +\delta  \\
 2 \alpha +2 \beta +\gamma +\delta  \\
 \alpha  \\
 \beta  \\
 -\alpha -2 \beta -\delta  \\
 -2 \alpha -\beta -\gamma  \\
 \gamma  \\
 \delta  \\
 \alpha +\gamma -\delta  \\
 \beta -\gamma +\delta  \\
 \beta -\gamma +\delta  \\
 \alpha +\gamma -\delta  \\
 -\alpha -\beta -\gamma -\delta  \\
 -2 \alpha -\beta -\gamma  \\
 -\alpha -2 \beta -\delta  \\
 -\alpha -\beta -\gamma -\delta  \\
 \delta  \\
 \gamma  \\
 \beta  \\
 \alpha 
\end{array}
\right).
\]

\bigskip

$X_3$ corresponds to the 5-dimensional embedding of dodecahedron with "broken faces". 
It would be interesting to understand the geometry of the $4$-dimensional embeddings corresponding to $X_4$ and $X_5.$ 
Since in the second case the opposite vertices identified it corresponds to the representation $V$ in agreement with (\ref{D}).

\bigskip


\newpage

{\bf The Truncated Cube}

\medskip

is one of the Archimedean solids, for which Buffon realisation is not affine equivalent to the standard one.

\begin{figure} [H]
   \centering
  \includegraphics[scale=0.7]{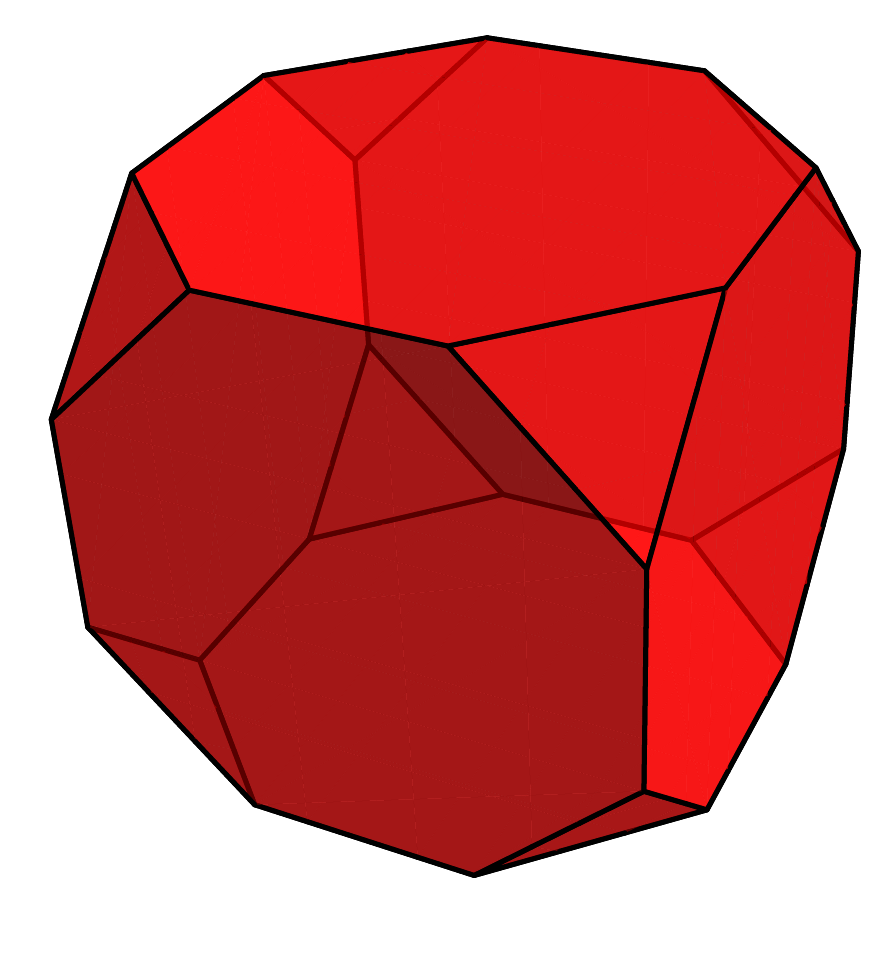}\quad{}\includegraphics[scale=1.0]{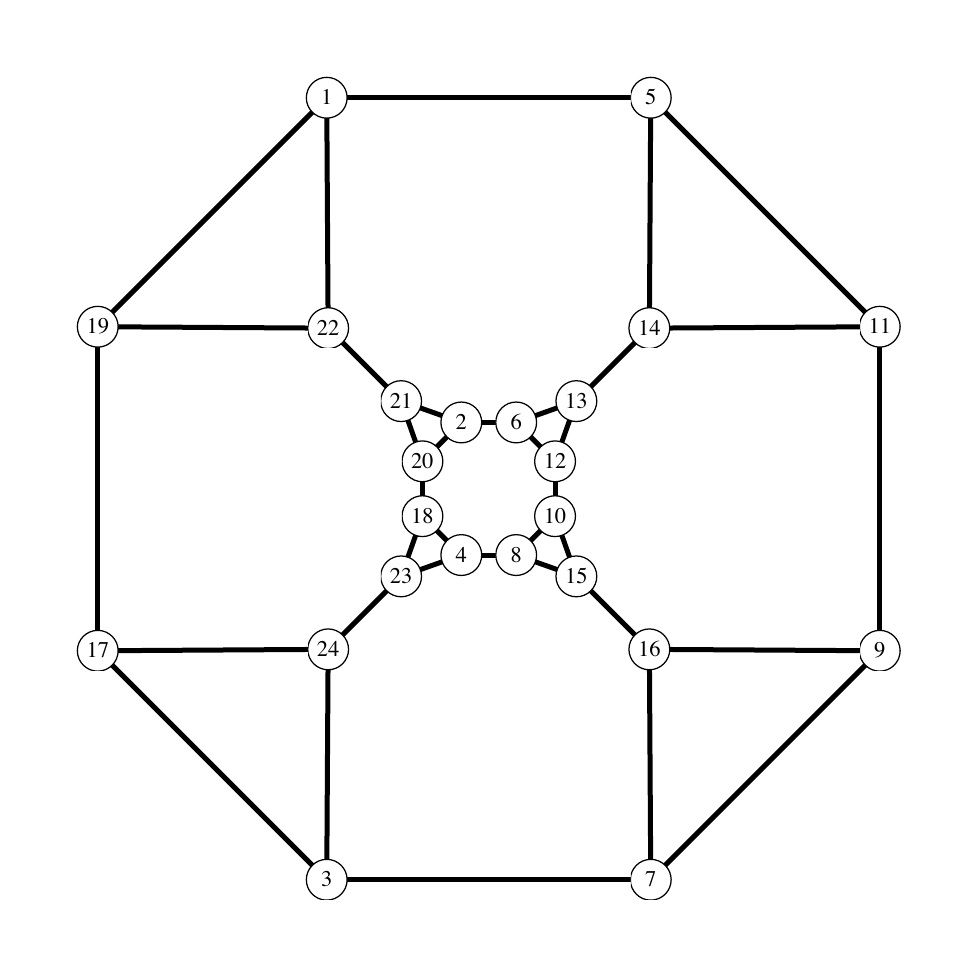}
  \bigskip
  \caption{The truncated cube and the corresponding $1$-skeleton graph.}
\end{figure}

The corresponding Buffon spectrum is:
\[
\left\{ 1^{(1)},\frac{1}{12}\left(7+\sqrt{17}\right)^{(3)},\frac{5}{6}^{(3)},\frac{2}{3}^{(1)},\frac{1}{2}^{(5)},\frac{1}{3}^{(3)},\frac{1}{12}\left(7-\sqrt{17}\right)^{(3)},\frac{1}{6}^{(5)}\right\} 
\]

The eigenspace corresponding to the second highest eigenvalue is:
\[
X_{2}=\left(\begin{array}{c}
\frac{1}{8}\left(\left(\sqrt{17}-1\right)(\alpha+2\gamma)+\left(\sqrt{17}-9\right)\beta\right)\\
\frac{1}{8}\left(\left(\sqrt{17}-1\right)(3\beta+2\gamma)-\left(7+\sqrt{17}\right)\alpha\right)\\
\frac{1}{8}\left(\left(7+\sqrt{17}\right)\alpha+\left(\sqrt{17}-9\right)\beta+2\left(\sqrt{17}-5\right)\gamma\right)\\
\frac{1}{8}\left(-\sqrt{17}\alpha+\alpha+3\left(\sqrt{17}-1\right)\beta+2\left(\sqrt{17}-5\right)\gamma\right)\\
\frac{1}{8}\left(\left(\sqrt{17}-1\right)\alpha-3\left(\sqrt{17}-1\right)\beta-2\left(\sqrt{17}-5\right)\gamma\right)\\
\frac{1}{8}\left(-\left(7+\sqrt{17}\right)\alpha-\left(\sqrt{17}-9\right)\beta-2\left(\sqrt{17}-5\right)\gamma\right)\\
\frac{1}{8}\left(\left(7+\sqrt{17}\right)\alpha-\left(\sqrt{17}-1\right)(3\beta+2\gamma)\right)\\
\frac{1}{8}\left((1-\sqrt{17})\alpha-\left(\sqrt{17}-9\right)\beta-2\left(\sqrt{17}-1\right)\gamma\right)\\
\frac{1}{8}\left(\left(\sqrt{17}-1\right)(3\alpha-2\gamma)-\left(7+\sqrt{17}\right)\beta\right)\\
\frac{1}{8}\left(\left(\sqrt{17}-9\right)\alpha+\left(\sqrt{17}-1\right)(\beta-2\gamma)\right)\\
\frac{1}{8}\left(-\left(\sqrt{17}-9\right)\alpha-\left(7+\sqrt{17}\right)\beta+2\left(\sqrt{17}-5\right)\gamma\right)\\
\frac{1}{8}\left(\left(\sqrt{17}-1\right)(-3\alpha+\beta)+2\left(\sqrt{17}-5\right)\gamma\right)\\
-\alpha\\
-\beta\\
-\gamma\\
\alpha-\beta-\gamma\\
\frac{1}{8}\left(3\left(\sqrt{17}-1\right)\alpha+(1-\sqrt{17})\beta-2\left(\sqrt{17}-5\right)\gamma\right)\\
\frac{1}{8}\left(\left(\sqrt{17}-9\right)\alpha+\left(7+\sqrt{17}\right)\beta-2\left(\sqrt{17}-5\right)\gamma\right)\\
\frac{1}{8}\left(-\left(\sqrt{17}-9\right)\alpha-\left(\sqrt{17}-1\right)(\beta-2\gamma)\right)\\
\frac{1}{8}\left(\left(\sqrt{17}-1\right)(-3\alpha+2\gamma)+\left(7+\sqrt{17}\right)\beta\right)\\
-\alpha+\beta+\gamma\\
\gamma\\
\beta\\
\alpha
\end{array}\right),
\]

The facing octagons of the geometric realisation of $X_{2}$ are not affine regular: one can check that 
$(x_{22}-x_{14})=\frac{3+\sqrt{17}}{4}(x_1-x_5)$ while for the regular octagon $(x_{22}-x_{14})=(1+\sqrt{2})(x_1-x_5).$
Thus the affine $B$-regular truncated cube obtained by the Buffon procedure is not an affine version of the regular truncated cube.  

\begin{figure} [H]
   \centering
  \includegraphics[scale=0.55]{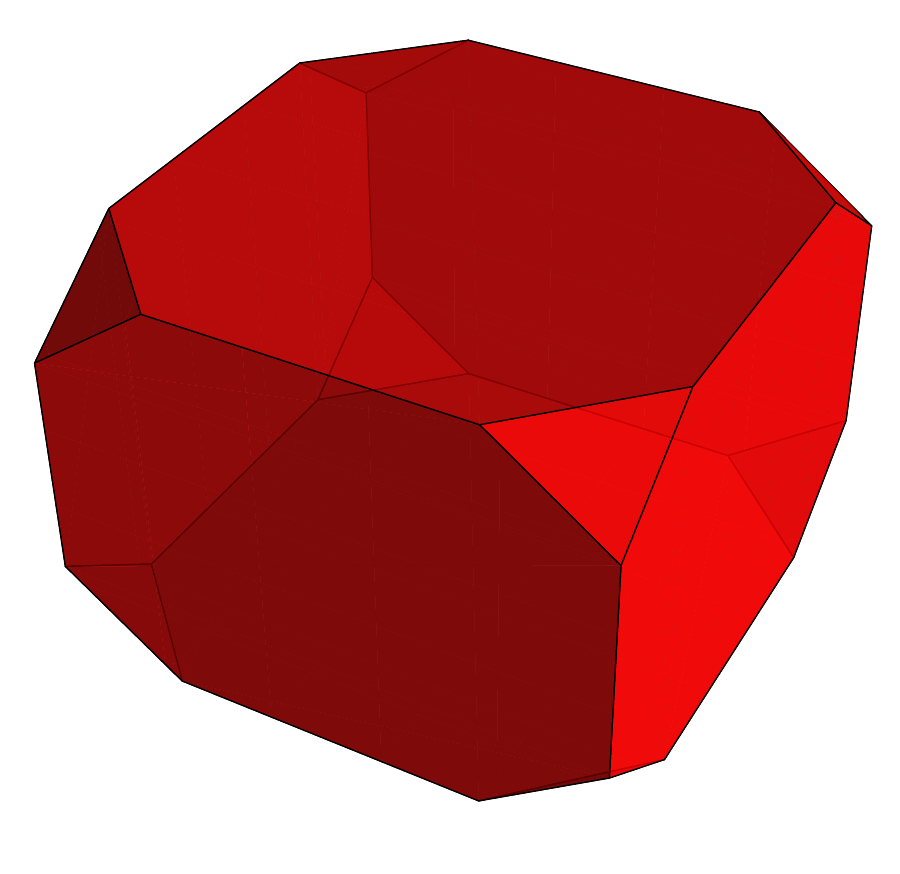}
  \caption{The affine $B$-regular truncated cube }
\end{figure}

\bigskip{}

\newpage

{\bf Triakis Tetrahedron}

\medskip

is the Catalan solids dual to the truncated tetrahedron. This is the simplest case 
when convexity does not hold for Buffon realisation.

\begin{figure} [H]
   \centering
  \includegraphics[scale=0.60]{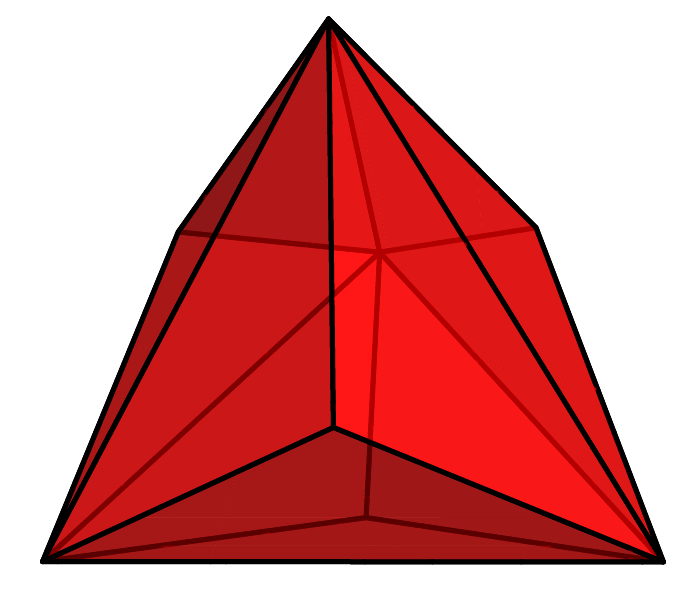}\quad{}\includegraphics[scale=0.55]{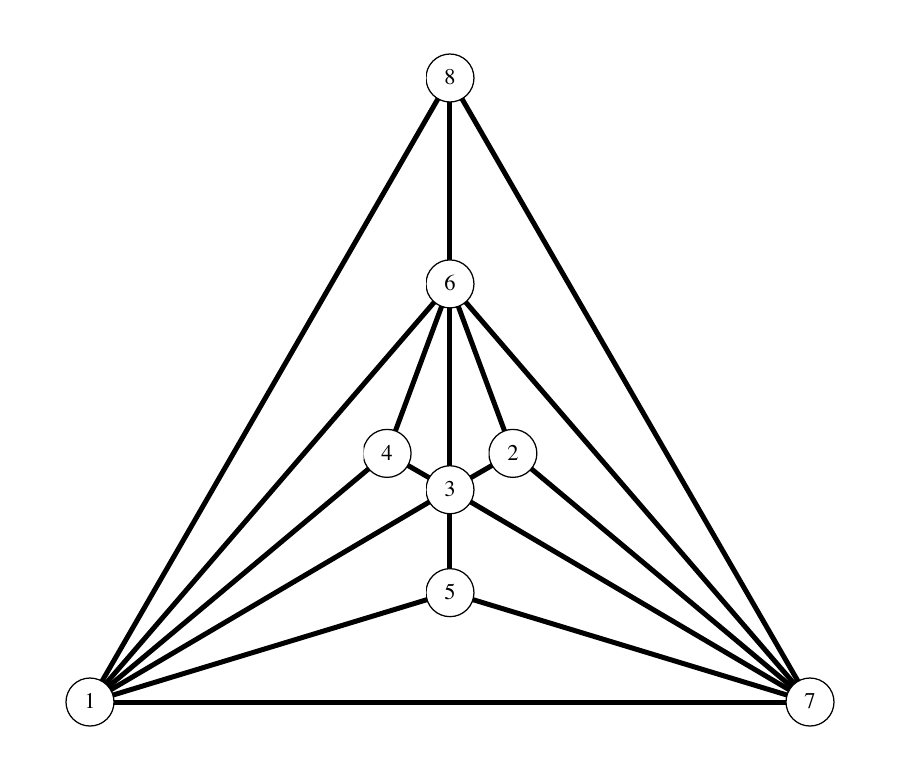}
\end{figure}

The corresponding Buffon eigenvalues are:
\[
\left\{ 1^{(1)},\frac{7}{12}^{(3)},\frac{1}{3}^{(3)},\frac{1}{4}^{(1)}\right\} 
\]
The eigenspaces corresponding to the eigenvalues $\lambda_{2}=\frac{7}{12}$
and to $\lambda_{3}=\frac{1}{3}$ are:

\[
X_{2}=\left(\begin{array}{c}
\frac{\alpha}{2}-\beta-\gamma\\
-\alpha+2\beta+2\gamma\\
-\frac{\alpha}{2}\\
-2\beta\\
-2\gamma\\
\gamma\\
\beta\\
\alpha
\end{array}\right)\quad X_{3}=\left(\begin{array}{c}
-\alpha-\beta-\gamma\\
-\alpha-\beta-\gamma\\
\alpha\\
\beta\\
\gamma\\
\gamma\\
\beta\\
\alpha
\end{array}\right)
\]

A particular geometric realisation derived from  $X_{2}$ is shown in Figure 12:

\begin{figure} [H]
   \centering
  \includegraphics[scale=0.6]{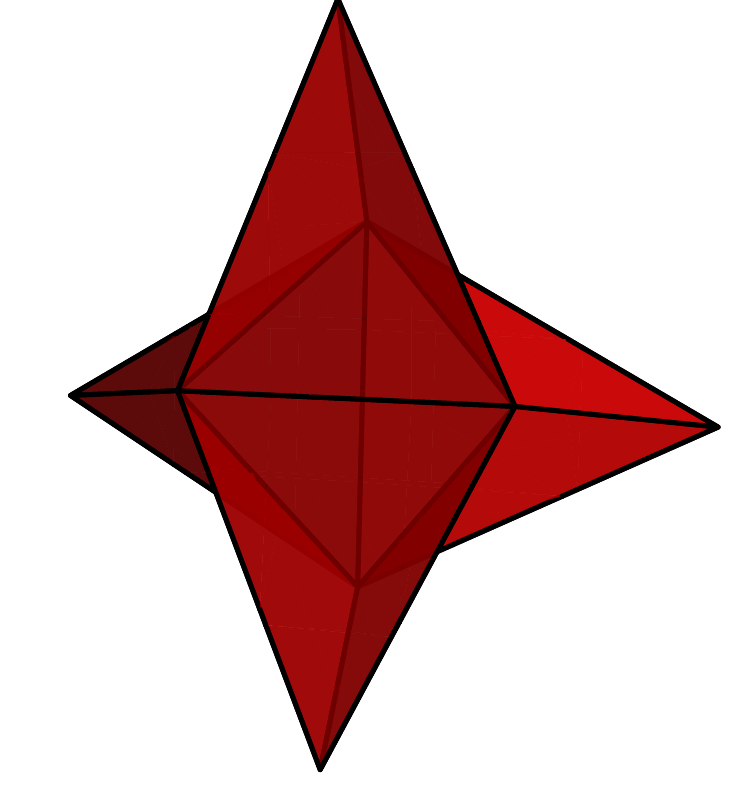}
   \caption{The affine $B$-regular version of the triakis tetrahedron is star-shaped but not convex.}
\end{figure}

As $X_{3}$ shows the vertices coalesce together pairwise. The corresponding geometrical realisation gives a general tetrahedron.

\medskip{}
\newpage

{\bf The Rhombic Dodecahedron}

\medskip

is the Catalan solid dual to cuboctahedron. We will see that it does not admit Buffon realisation.

\begin{figure} [H]
   \centering
  \includegraphics[scale=0.6]{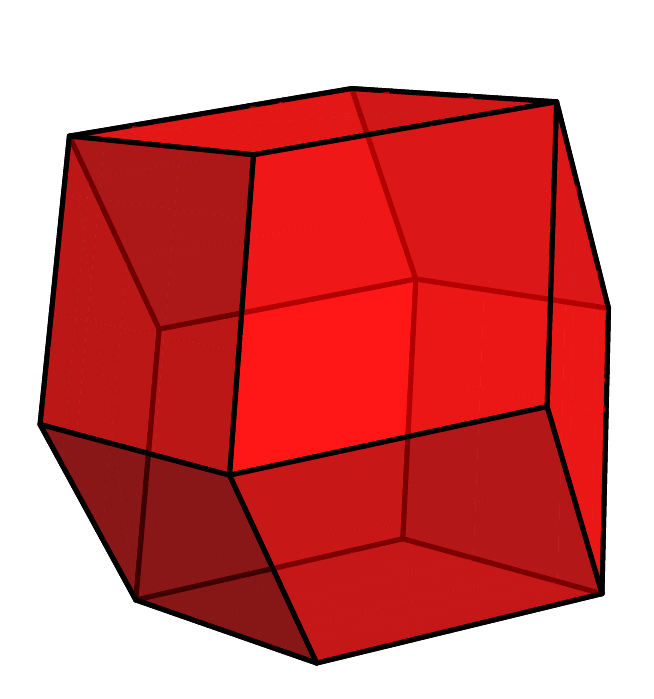}\quad{}\includegraphics[scale=0.55]{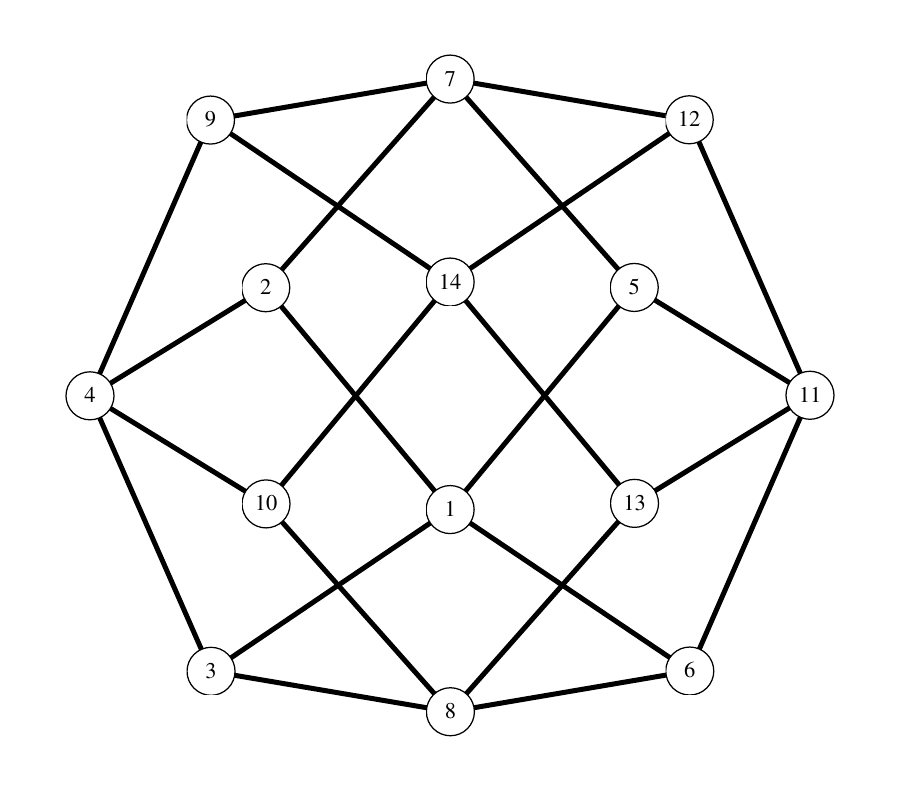}
   \caption{The rhombic dodecahedron and the corresponding $1$-skeleton graph.}
\end{figure}

The corresponding eigenvalues are:
\[
\left\{ 1^{(1)},\frac{1}{6}\left(3+\sqrt{3}\right)^{(3)},\frac{1}{2}^{(6)},\frac{1}{6}\left(3-\sqrt{3}\right)^{(3)},0^{(1)}\right\} 
\]
The $3$-dimensional eigenspaces corresponding to the eigenvalues
$\lambda_{2}=\frac{1}{6}\left(3+\sqrt{3}\right)$ and $\lambda_{4}=\frac{1}{6}\left(3-\sqrt{3}\right)$
are:

\[
X_{2}=\left(\begin{array}{c}
-\alpha\\
-\beta\\
-\gamma\\
\alpha-\frac{\sqrt{3}\beta}{2}-\frac{\sqrt{3}\gamma}{2}\\
\gamma-\frac{2\alpha}{\sqrt{3}}\\
\beta-\frac{2\alpha}{\sqrt{3}}\\
\frac{\sqrt{3}\gamma}{2}-\frac{\sqrt{3}\beta}{2}\\
\frac{\sqrt{3}\beta}{2}-\frac{\sqrt{3}\gamma}{2}\\
\frac{2\alpha}{\sqrt{3}}-\beta\\
\frac{2\alpha}{\sqrt{3}}-\gamma\\
-\alpha+\frac{\sqrt{3}\beta}{2}+\frac{\sqrt{3}\gamma}{2}\\
\gamma\\
\beta\\
\alpha
\end{array}\right)\quad X_{4}=\left(\begin{array}{c}
-\alpha\\
-\beta\\
-\gamma\\
\alpha+\frac{\sqrt{3}\beta}{2}+\frac{\sqrt{3}\gamma}{2}\\
\frac{2\alpha}{\sqrt{3}}+\gamma\\
\frac{2\alpha}{\sqrt{3}}+\beta\\
\frac{\sqrt{3}\beta}{2}-\frac{\sqrt{3}\gamma}{2}\\
\frac{\sqrt{3}\gamma}{2}-\frac{\sqrt{3}\beta}{2}\\
-\frac{2\alpha}{\sqrt{3}}-\beta\\
-\frac{2\alpha}{\sqrt{3}}-\gamma\\
-\alpha-\frac{\sqrt{3}\beta}{2}-\frac{\sqrt{3}\gamma}{2}\\
\gamma\\
\beta\\
\alpha
\end{array}\right)
\]

Both eigenspaces $X_{2}$ and  $X_{4}$ fail to correspond to a polyhedron with the combinatorial structure of the $1$-skeleton of the rhombic dodecahedron. A particular graph realisation obtained from the eigenspace  $X_{2}$ is shown in Fig. 13.

\begin{figure} [H]
   \centering
  \includegraphics[scale=0.6]{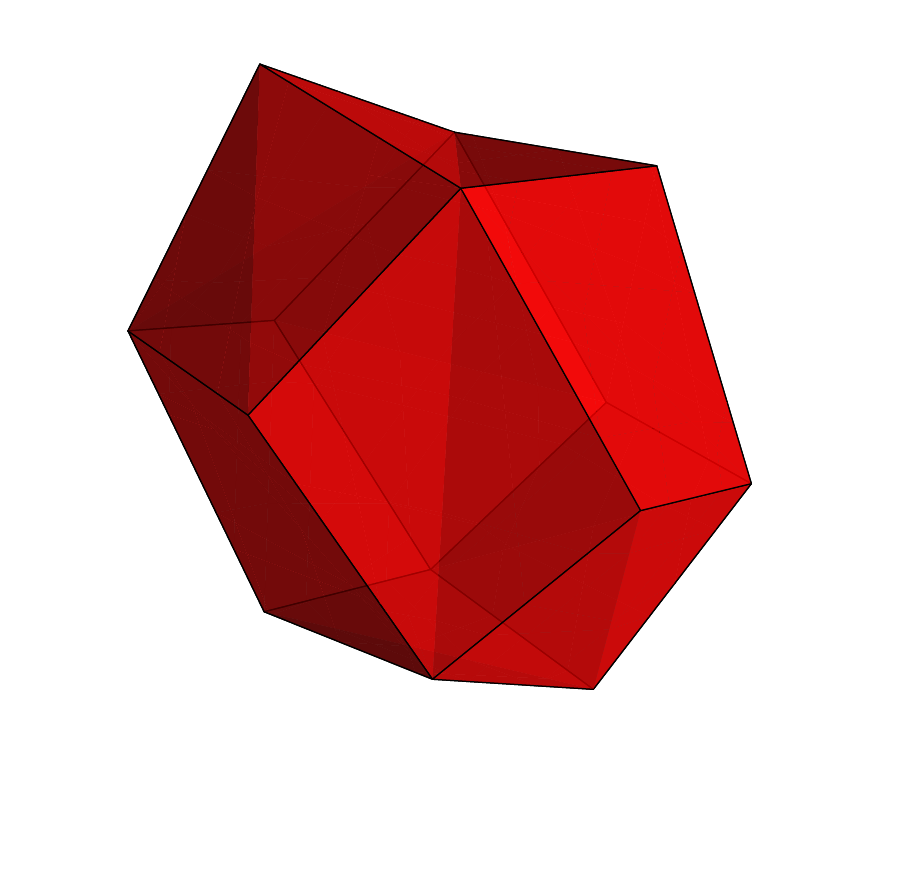}
  \caption{Subdominant eigenspace realisation of the rhombic dodecahedron. All faces are broken (non-planar). }
\end{figure}

{\bf Pentakis Dodecahedron}

\medskip

is the Catalan solid dual to the truncated icosahedron, which we mentioned in the Introduction.

\begin{figure} [H]
   \centering
  \includegraphics[scale=0.7]{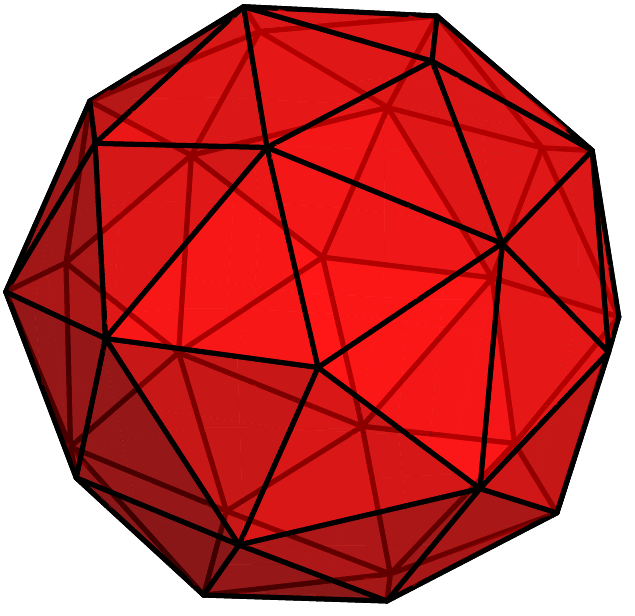}\includegraphics[scale=0.6]{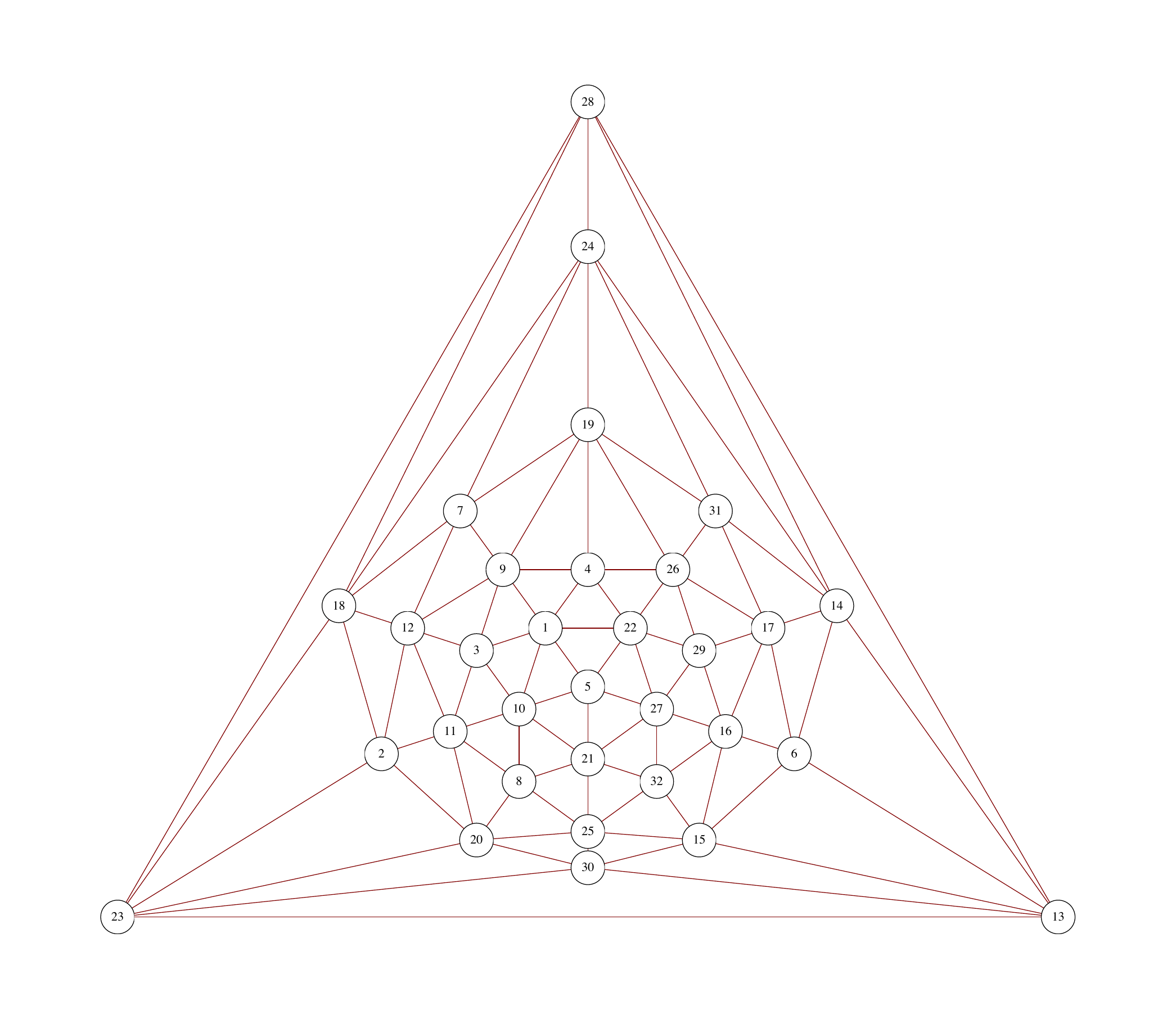}
  \medskip
  \caption{The pentakis dodecahedron and the corresponding $1$-skeleton graph.}
\end{figure}

The corresponding Buffon eigenvalues are:
$$
1^{(1)},
\frac{1}{120} \left(60+5 \sqrt{5}+\sqrt{725+240 \sqrt{5}}\right)^{(3)},
\frac{1}{120} \left(65+\sqrt{385}\right)^{(5)},
\frac{1}{24}   \left(12-\sqrt{5}+\sqrt{29-\frac{48}{\sqrt{5}}}\right)^{(3)}, 
\frac{1}{2}^{(4)}, 
$$
$$
\frac{1}{120} \left(65-\sqrt{385}\right)^{(5)},
\frac{1}{3}^{(4)},
\frac{1}{120} \left(60+5 \sqrt{5}-\sqrt{725+240 \sqrt{5}}\right)^{(3)},
\frac{1}{24}   \left(12-\sqrt{5}-\sqrt{29-\frac{48}{\sqrt{5}}}\right)^{(3)},
\frac{1}{4}^{(1)}.
$$

The subdominant eigenspace corresponds to:
$$
\left(
\begin{array}{c}
 \phi \left(\left(\sqrt{5}-3\right) \beta +(1-\sqrt{5}) \gamma \right) \\
 -\alpha -\frac{1}{2} \left(\sqrt{5}-1\right) (\beta -\gamma ) \\
 \frac{1}{2} \left((1-\sqrt{5}) \alpha -2 \beta +(1-\sqrt{5}) \gamma \right) \\
 -\gamma  \\
 \frac{1}{2} \left(\left(\sqrt{5}-1\right) \alpha +(1-\sqrt{5}) \beta -2 \gamma \right) \\
 \frac{1}{2} \left(\left(\sqrt{5}-1\right) \alpha +2 \beta +\left(\sqrt{5}-1\right) \gamma \right) \\
 -\alpha  \\
 -\beta  \\
 \phi \left(\left(\sqrt{5}-3\right) \alpha -2
   \left(\sqrt{5}-2\right) \beta +\left(\sqrt{5}-3\right) \gamma \right) \\
 - \phi \left(\left(\sqrt{5}-1\right) \beta
   -\left(\sqrt{5}-3\right) \gamma \right) \\
 \phi \left(\left(\sqrt{5}-3\right) \alpha +(1-\sqrt{5}) \beta
  \right) \\
 - \phi \left(\left(\sqrt{5}-1\right) \alpha
   -\left(\sqrt{5}-3\right) \beta \right) \\
 - \phi \left(\left(\sqrt{5}-3\right) \beta +(1-\sqrt{5}) \gamma
   \right) \\
  \phi \left(\left(\sqrt{5}-1\right) \beta
   -\left(\sqrt{5}-3\right) \gamma \right) \\
 - \phi \left(\left(\sqrt{5}-3\right) \alpha -2
   \left(\sqrt{5}-2\right) \beta +\left(\sqrt{5}-3\right) \gamma \right) \\
  \phi \left(\left(\sqrt{5}-1\right) \alpha
   -\left(\sqrt{5}-3\right) \beta \right) \\
 - \phi \left(\left(\sqrt{5}-3\right) \alpha +(1-\sqrt{5}) \beta
   \right) \\
 - \phi \left(\left(\sqrt{5}-1\right) \alpha
   +\left(\sqrt{5}-3\right) \gamma \right) \\
-2 \phi \left(\sqrt{5}-2\right) (\alpha -\beta +\gamma ) \\
 - \phi \left(2 \left(\sqrt{5}-2\right) \alpha
   -\left(\sqrt{5}-3\right) (\beta -\gamma )\right) \\
 - \phi \left(\left(\sqrt{5}-3\right) \alpha
   -\left(\sqrt{5}-3\right) \beta +2 \left(\sqrt{5}-2\right) \gamma \right) \\
 - \phi \left(\left(\sqrt{5}-3\right) \alpha
   +\left(\sqrt{5}-1\right) \gamma \right) \\
  \phi \left(\left(\sqrt{5}-3\right) \alpha
   +\left(\sqrt{5}-1\right) \gamma \right) \\
 \phi \left(\left(\sqrt{5}-3\right) \alpha
   -\left(\sqrt{5}-3\right) \beta +2 \left(\sqrt{5}-2\right) \gamma \right) \\
 2 \phi \left(\sqrt{5}-2\right)  (\alpha -\beta +\gamma ) \\
 \phi   \left(2 \left(\sqrt{5}-2\right) \alpha
   -\left(\sqrt{5}-3\right) (\beta -\gamma )\right) \\
 \phi  \left(\left(\sqrt{5}-1\right) \alpha
   +\left(\sqrt{5}-3\right) \gamma \right) \\
 \frac{1}{2} \left((1-\sqrt{5}) \alpha +\left(\sqrt{5}-1\right) \beta +2 \gamma \right) \\
 \alpha +\frac{1}{2} \left(\sqrt{5}-1\right) (\beta -\gamma ) \\
 \gamma  \\
 \beta  \\
 \alpha 
\end{array}
\right)
$$
where $\phi= \frac{1}{24} (5+\sqrt{145+48 \sqrt{5}}).$

The corresponding Buffon realization is shown below. It is convex and looks pretty similar to the Catalan version, but the pyramids are slightly higher. The ratios of the height of a pyramid to the distance of its top vertex from the centre in the Catalan and Buffon cases are $1/3 (1- 1/\sqrt{5})\approx 0.184$ and $1-1/12 (\sqrt{5} + \sqrt{29 + 48/\sqrt{5}}) \approx 0.222$ respectively.

\begin{figure} [H]
   \centering
  \includegraphics[scale=0.9]{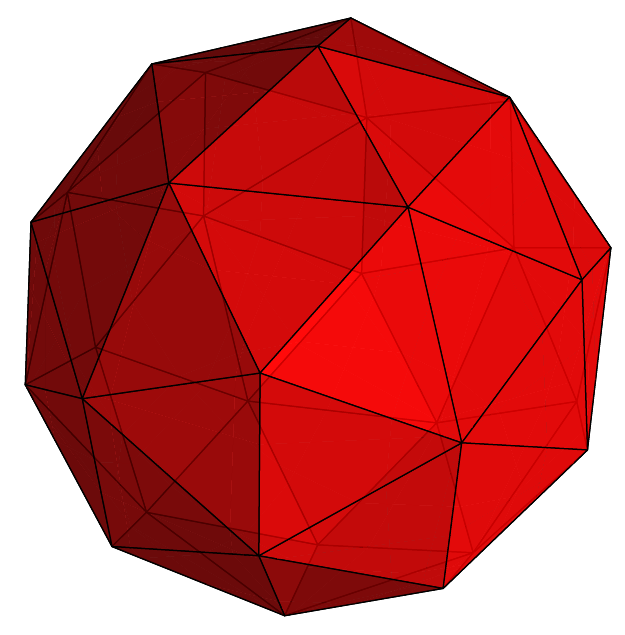}
  \caption{Affine $B$-regular pentakis dodecahedron.}
\end{figure}

The self-intersecting realisations corresponding to other multiplicity 3 eigenvalues are shown below.

\begin{figure}[H]
\centerline{\includegraphics[scale=1.2]{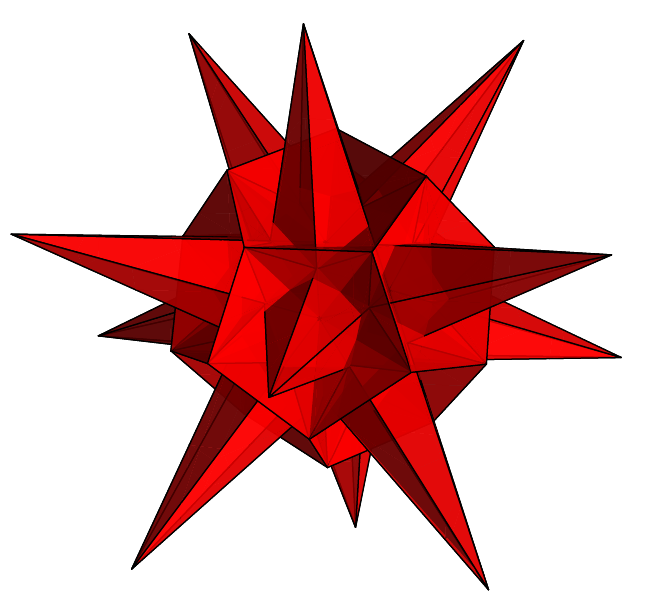}}
\caption{ Eigenspace realisation corresponding to $\lambda_8=(60+5 \sqrt{5}-\sqrt{725+240 \sqrt{5}})/120$ conjugated to $\lambda_2$.
The pyramids are build inside and go through the dodecahedron.} 
\end{figure}

\begin{figure}[H]
\centerline{ \includegraphics[scale=0.35]{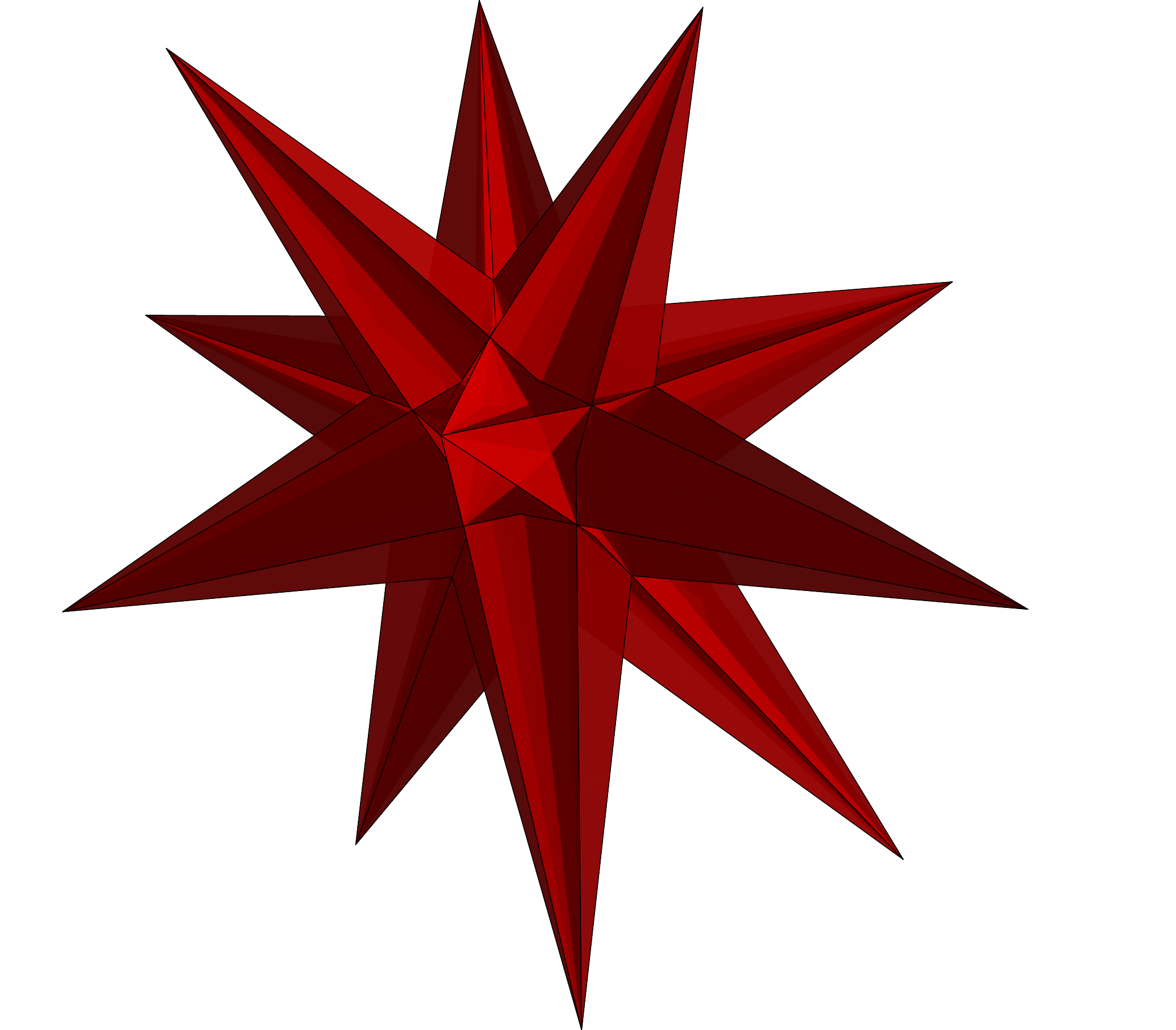} \quad \includegraphics[scale=0.3]{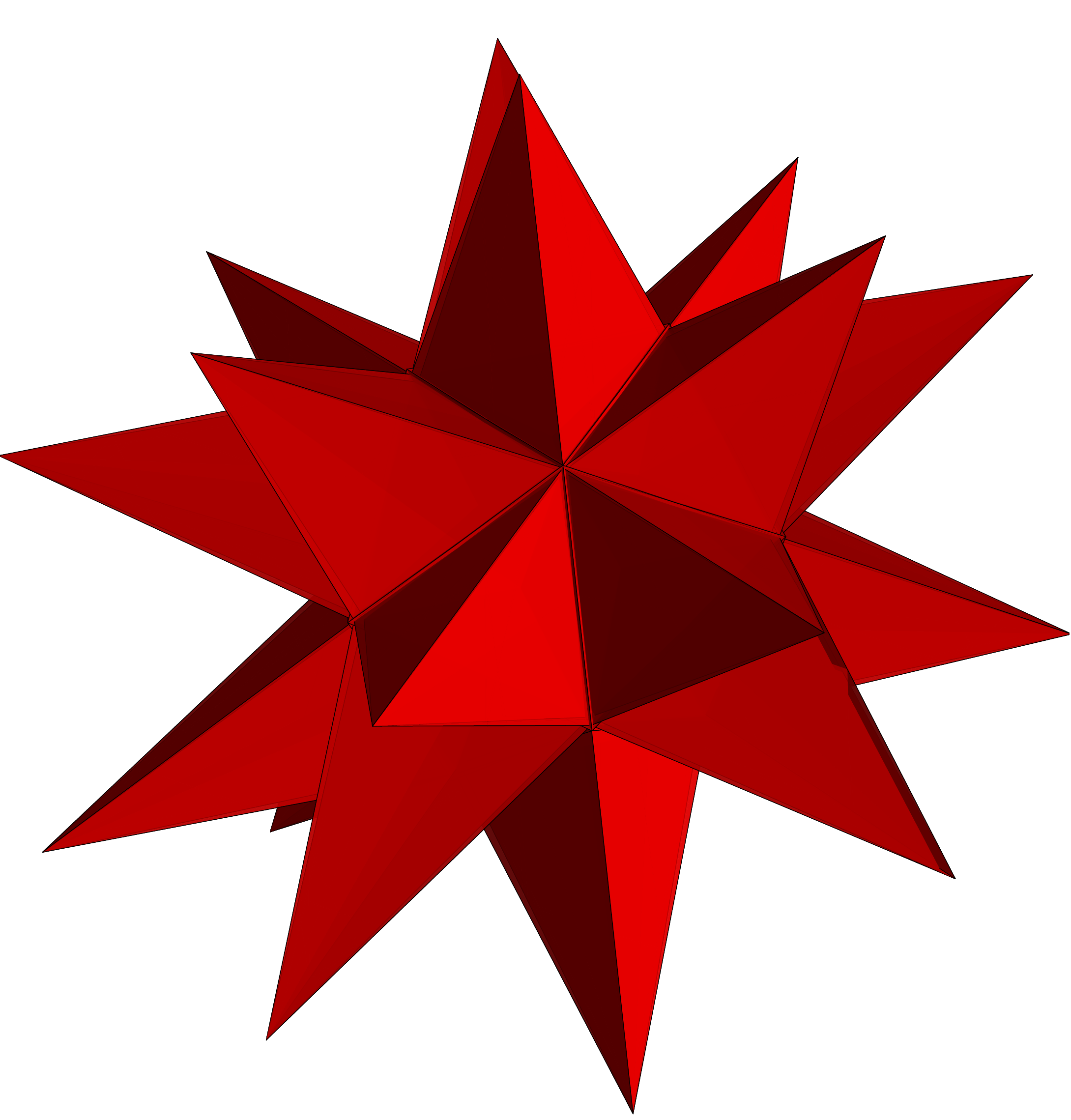} }
\caption{ Eigenspace realisations corresponding to two remaining multiplicity 3 eigenvalues $\lambda_4$ and $\lambda_9$: great icosahedron and great stellated dodecahedron with extra vertices.} 
\end{figure}

To finish we present Mathematica image of the cumulated dodecahedron with all edges of equal length, featured on Leonardo's drawing (Fig.2):

\begin{figure}[H]
\centerline{\includegraphics[scale=0.8]{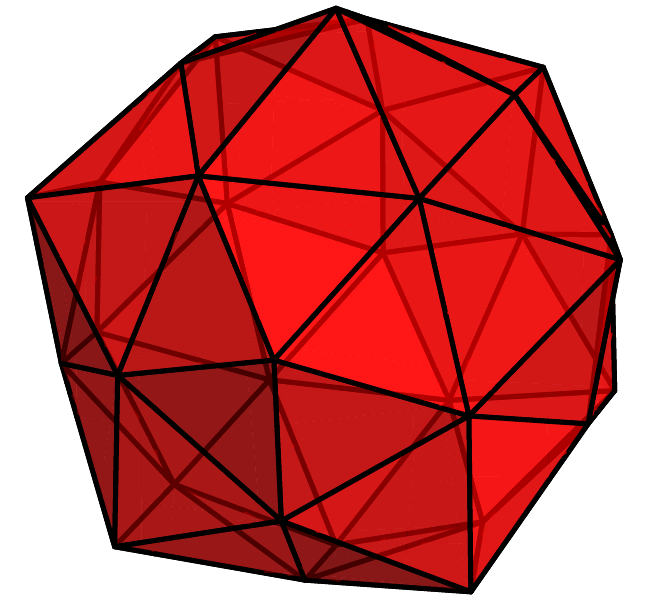}}
\caption{Cumulated dodecahedron.} 
\end{figure}
\end{document}